\documentclass[11pt,reqno,a4paper]{article}

\usepackage{geometry}
\usepackage{amsmath,amsthm,amstext,amscd,amssymb,euscript,url}
\usepackage{mathrsfs}
\usepackage{calrsfs}
\usepackage{epsfig}
\usepackage[inline]{enumitem}
\usepackage{mathtools}
\usepackage[absolute,overlay]{textpos}
\usepackage{graphicx}
\usepackage{subfig}
\usepackage{tikz}
\usetikzlibrary{matrix}
\usetikzlibrary{arrows,positioning}

\usepackage{color}

\newcommand{\Z}{\mathbb Z}
\newcommand{\R}{\mathbb R}
\newcommand{\N}{\mathbb N}

\newcommand{\E}{\mathbb E}
\newcommand{\Zd}{\mathbb Z^d}

\renewcommand{\phi}{\varphi}

\newcommand{\loc}{\mathcal{L}}

\def\1{{\mathchoice {\rm 1\mskip-4mu l} {\rm 1\mskip-4mu l}
{\rm 1\mskip-4.5mu l} {\rm 1\mskip-5mu l}}}

\newtheorem{theorem}{{\small T}{\scriptsize HEOREM}}[section]
\newtheorem{corollary}{{\bf{\small C}{\scriptsize OROLLARY}}}[section]
\newtheorem{proposition}{{\bf{\small P}{\scriptsize ROPOSITION}}}[section]
\newtheorem{lemma}{{\bf{\small L}{\scriptsize EMMA}}}[section]
\newtheorem{remark}{{\bf{\small R}{\scriptsize EMARK}}}[section]
\newtheorem{definition}{{\bf{\small D}{\scriptsize EFINITION}}}[section]

\renewenvironment{proof}[1]
{\noindent{{\bf{\small{ P}{\scriptsize ROOF}}}.}\hspace{0.1cm} #1} {$\;\qed$\newline}

\newcommand{\beq}{\begin{eqnarray}}
\newcommand{\eeq}{\end{eqnarray}}

\newcommand{\ba}{\begin{align*}}
\newcommand{\ea}{\end{align*}}

\newcommand{\be}{\begin{equation}}
\newcommand{\ee}{\end{equation}}

\newcommand{\bl}{\begin{lemma}}
\newcommand{\el}{\end{lemma}}

\newcommand{\br}{\begin{remark}}
\newcommand{\er}{\end{remark}}

\newcommand{\bt}{\begin{theorem}}
\newcommand{\et}{\end{theorem}}

\newcommand{\bd}{\begin{definition}}
\newcommand{\ed}{\end{definition}}

\newcommand{\bp}{\begin{proposition}}
\newcommand{\ep}{\end{proposition}}

\newcommand{\bc}{\begin{corollary}}
\newcommand{\ec}{\end{corollary}}

\newcommand{\bpr}{\begin{proof}}
\newcommand{\epr}{\end{proof}}

\newcommand{\bi}{\begin{itemize}}
\newcommand{\ei}{\end{itemize}}

\newcommand{\ben}{\begin{enumerate}}
\newcommand{\een}{\end{enumerate}}

\newcommand{\caA}{{\mathcal A}}
\newcommand{\caB}{{\mathcal B}}

\newcommand{\caE}{{\mathrsfs E}}

\newcommand{\caF}{{\mathcal F}}

\newcommand{\caR}{{\mathcal R}}

\newcommand{\caX}{{\mathcal X}}

\newcommand{\caZ}{{\mathscr Z}}

\renewcommand{\(}{\left(}
\renewcommand{\)}{\right)}
\newcommand{\nn}{\nonumber}

\newmuskip\pFqmuskip

\newcommand\pFq[6][8]{%
	\begingroup 
	\pFqmuskip=#1mu\relax
	\mathcode`\,=\string"8000
	\begingroup\lccode`\~=`\,
	\lowercase{\endgroup\let~}\pFqcomma
	{}_{#2}F_{#3}{\left[\genfrac..{0pt}{}{#4}{#5};#6\right]}%
	\endgroup
}
\newcommand{\pFqcomma}{\mskip\pFqmuskip}

\newcommand*{\myprime}{^{\prime}\mkern-1.2mu}
\newcommand*{\mydprime}{^{\prime\prime}\mkern-1.2mu}

\newcommand{\quotes}[1]{``#1''}
\newcommand{\norm}[1]{\left\lVert#1\right\rVert}

\newtheorem{asu}{Assumption}
\newcommand{\basu}{\begin{asu}}
\newcommand{\easu}{\end{asu}}

\DeclareMathOperator\erf{erf}

\begin{document}
\title{{\bf Condensation of SIP particles and sticky Brownian motion}}
\author{Mario Ayala, Gioia Carinci and Frank Redig
\\
\small{Delft Institute of Applied Mathematics}\\
\small{Delft University of Technology}\\
{\small Van Mourik Broekmanweg 6, 2628 XE Delft}
\\
\small{The Netherlands}
}
\maketitle

\begin{abstract}
We study the symmetric inclusion process (SIP) in the condensation regime. We obtain an explicit scaling  for the variance of the density field in this regime, when initially started from a homogeneous product measure. This provides relevant new information on the coarsening dynamics of condensing interacting particle systems on the infinite lattice.\\
We obtain our result by proving convergence to sticky Brownian motion for the difference of positions of two SIP particles in the sense of Mosco convergence of Dirichlet forms. Our approach implies the convergence of the probabilities of two SIP particles to be together at time $t$. This, combined with self-duality, allows us to obtain the explicit scaling for the variance of the fluctuation field. 
\end{abstract}

\tableofcontents

\section{Introduction}

The symmetric inclusion process (SIP) is an interacting particle system where a single particle performs symmetric continuous-time random walks on the lattice $\Z$, with rates $k p(i,j)=k p(j,i)$ ($k>0$) and where particles interact by attracting each other (see below for the precise definition) at rate $p(i,j) \eta_i \eta_j$, where $\eta_i$ is the number of particles at site $i$. The parameter $k$ regulates the relative strength of diffusion w.r.t. attraction between particles. The symmetric inclusion process is self-dual, and many results on its macroscopic behavior can be obtained via this property. Self-duality implies that the expectation of the number of particles can be understood  from one dual particle. In particular, because one dual particle scales to Brownian motion in the diffusive scaling, the hydrodynamic limit of SIP is the heat equation. The next step is to understand the variance of the density field, which requires two dual particles.

It is well-known that in the regime $k\to 0$ the SIP manifests condensation (the attractive interaction dominates), and via the  self-duality of  SIP more information can be obtained about this condensation process than for a generic process (such as zero range processes). Indeed, in  \cite{carinci2017exact} two of the authors of this paper in collaboration with C. Giardin\`{a} have obtained an explicit formula for the Fourier-Laplace transform of two particle transition probabilities for interacting particle systems such as the simple symmetric exclusion and the simple symmetric inclusion process, where simple refers to nearest neighbor in dimension 1. From this formula, the authors were able to extract information about the variance of the time-dependent density field in starting from a homogeneous product measure. With the help of duality this reduces to the study of the scaling behavior of two dual particles. In particular, for the inclusion process in the condensation regime,  from the study of the scaling behavior of the time dependent variance of the density field, one can extract information about the coarsening process. It turned out that the scaling limit of two particles is in that case a pair of sticky Brownian motions. From this one can infer the qualitative picture that in the condensation regime, when started from a homogeneous product measure, large piles of particles are formed which move as Brownian motion, and interact with each other as sticky Brownian motions.

The whole analysis in \cite{carinci2017exact} is based on the exact formula for the Fourier-Laplace transform of the transition probabilities of two SIP particles as mentioned above. This exact computation is based on the fact that the underlying random walk is nearest neighbor, and therefore the results are restricted to that case. However, we expect that for the SIP in the condensation regime, sticky Brownian motion appears as a scaling limit in much larger generality in dimension 1. The exact formula in \cite{carinci2017exact} yields convergence of semigroups, and therefore convergence of finite dimensional distributions. However, because of the rescaling in the condensation regime, one cannot expect convergence of generators, but rather a convergence result in the spirit of slow-fast systems, i.e., of the type gamma convergence. Moreover, the difference of two SIP-particles is not simply a random walk slowed down when it is the origin as in e.g. \cite{amir1991sticky}. Instead, it is a random walk which is pulled towards the origin when it is close to it, which only in the scaling limit leads to a slow-down at the origin, i.e., sticky Brownian motion.

In this paper,  we obtain a precise scaling behavior of the variance of the density field in the condensation regime. We find the explicit scaling form for this variance in real time (as opposed to the Laplace transformed result in \cite{carinci2017exact}), thus giving more insight in the coarsening process when initially started from a homogeneous product measure of density $\rho$. This is the first rigorous result on coarsening dynamics in interacting particle systems directly on infinite lattices, for a general class of underlying random walks. There exist important results on condensation either heuristically on the infinite lattice or rigorous but constrained to finite lattices. For example \cite{cao2014dynamics} heuristically  discusses on infinite lattices the effective motion of clusters in the coarsening process for the TASIP; or the work \cite{chau2015explosive} which based on heuristic mean field arguments studies the coarsening regime for the explosive condensation model. On the other hand, on finite lattices  via martingale techniques \cite{beltran2017martingale} studies the evolution of a condensing zero range process. In the context of the SIP, the authors of \cite{grosskinsky2013dynamics} on a finite lattice, showed the emergence of condensates as the parameter $k \to 0$ and rigorously characterize their dynamics. We also mention the recent work \cite{jatuviriyapornchai2020structure} where the structure of the condensed phase in SIP is analyzed in stationarity, in the thermodynamic limit.

Our main result is obtained by proving that the difference of two SIP particles converges to a two-sided sticky Brownian motion in the sense of Mosco convergence of Dirichlet forms originally introduced in \cite{mosco1994composite} and extended to the case of varying state spaces in \cite{kuwae2003convergence}. Because this notion of convergence implies convergence of semigroups in the $L^2$ space of the reversible measure, which is $dx + \gamma \delta_0$ for the sticky Brownian motion with stickiness parameter $\gamma>0$, the convergence of semigroups also implies that of transition probabilities of the form $p_t(x,0)$. This, together with self-duality, helps to explicitly obtain the limiting variance of the fluctuation field. Technically speaking, the main difficulty in our approach is that we have to define carefully how to transform functions defined on the discretized rescaled lattices into functions on the continuous limit space in order to obtain convergence of the relevant Hilbert spaces, and at the same time obtain the second condition of Mosco convergence. Mosco convergence is a weak form of convergence which is not frequently used in the probabilistic context. In our context it is however exactly the form of convergence which we need to study the variance of the density field. As already mentioned before, as it is strongly related to gamma-convergence, it is also a natural form of convergence in a setting reminiscent of slow-fast systems.

The rest of our paper is organized as follows. In Section 2 we deal with some preliminary notions; we introduce both the inclusion and the difference process in terms of their infinitesimal generators. In this section we also introduce the concept of duality and describe the appropriate regime in which condensation manifests itself. Our main result is stated in Section 3, were we present some non-trivial information about the variance of the time-dependent density field in the condensation regime and provide some heuristics for the dynamics described by this result.  Section 4 deals with the basic notions of Dirichlet forms. In the same section, we also introduce the notion of Mosco convergence on varying Hilbert spaces together with some useful simplifications in our setting. In Section 5, we present the proof of our main result and also show that the finite range difference process converges in the sense of Mosco convergence of Dirichlet forms to the two sided sticky Brownian motion. Finally, as supplementary material in the Appendix, we construct via stochastic time changes of Dirichlet forms the two sided sticky Brownian motion at zero and we also deal with the convergence of independent random walkers to standard Brownian motion. This last result, despite of being basic becomes a corner stone for our results of Section 5.

\section{Preliminaries}

\subsection{The Model: inclusion process}

The Symmetric Inclusion Process (SIP) is an interacting particle system  where particles randomly hop on the lattice $\Z$ with attractive interaction and no restrictions on the number of particles per site.
Configurations are denoted by $\eta$  and are elements of $\Omega=\N^{\Z}$ (where
$\N$ denotes the set of natural numbers including zero).
We denote by $\eta_x$ the number of particles at position $x \in \Z$  in the configuration $\eta\in\Omega$.
The generator working on local functions $f:\Omega\to\R$ is of the type
\be\label{SIPgen}
\loc f(\eta) = \sum_{i,j \in \Z} p(j-i) \eta_i (k +\eta_j) (f(\eta^{ij})-f(\eta))
\ee
where $\eta^{ij}$ denotes the configuration obtained from $\eta$ by removing a particle from $i$ and putting it at $j$. For the associated Markov process on $\Omega$, we use the notation $\{\eta(t):t\geq 0\}$, i.e., $\eta_x(t)$ denotes the number of particles at time $t$ at location $x \in \Z$. Additionally, we assume that the function $p:\R\to[0,\infty)$ satisfies the following properties
\ben
\item Symmetry: $p(r)=p(-r)$ for all $r\in \R$
\item Finite range: there exists $R>0$ such that: $p(r)=0$ for all $|r|>R$.
\item Irreducibility: for all $x,y\in\Z$ there exists  $n\in \N$ and  $x=i_1, i_2, \ldots,i_{n-1},i_n=y$, such that
$\prod\limits_{k=1}^{n-1} p(i_{k+1}-i_k)>0$.
\een
It is known that these particle systems have a one parameter family of homogeneous (w.r.t. translations) reversible and ergodic product measures $\mu_{\rho}, \rho>0$ with marginals
\[
\mu_\rho (\eta_i = n)= \frac{k^k \rho^n }{(k +\rho)^{k+n}} \frac{ \Gamma(k+n)}{\Gamma(n+1)\Gamma(k)}
\]
This family of measures is indexed by the density of particles, i.e.,
\[
\int\eta_0 d \mu_{\rho}= \rho
\]

\br
Notice that for these systems the initial configuration has to be chosen in a subset of configurations such that the process $\{\eta(t):t\geq 0\}$ is well-defined. A possible such subset is the set of tempered configurations. This is the set of configurations $\eta$ such that there exist $C, \beta \in \R$ that satisfy $ |\eta (x)| \leq C |x|^\beta$ for all $x\in \R$. 
We denote this set (with slight abuse of notation) still by $\Omega$, because
we will always start the process from such configurations, and this set has $\mu_{\bar{\rho}}$ measure $1$ for all $\rho$. Since
we are working mostly in $L^2(\mu_{\bar{\rho}})$ spaces, this is not a restriction.
\er

\subsection{Self-duality}

Let us denote by $\Omega_f \subseteq \Omega$ the set of configurations with a finite number of particles. We then have the following definition:

\bd
 We say that the process $\{ \eta_t : t \geq 0  \}$ is self-dual with self-duality function $D:\Omega_f\times \Omega \to\R$ if 
\be\label{dual1}
\E_\eta \big[D(\xi,\eta_t)\big]=\E_\xi \big[D(\xi_t, \eta)\big]
\ee
for all $ t \geq 0$ and $\xi\in \Omega_f, \eta\in \Omega$.
\ed
\noindent
In the definition above $\E_\eta$ and $\E_\xi$ denote expectation when the processes $\{ \eta_t : t \geq 0  \}$ and $\{ \eta_t : t \geq 0  \}$ are initialized from the configuration $\eta$ and $\xi$ respectively . Additionally we require the duality functions to be of factorized form, i.e.,
\be\label{dual2}
D(\xi,\eta)=\prod_{i \in \Z} d(\xi_i,\eta_i)
\ee
where the single site duality function $d(m,\cdot)$ is a polynomial of degree $m$, more precisely
\be\label{singledual}
d(m,n)= \frac{n! \Gamma(k)}{(n-m)! \Gamma(k+m) } \1_{\{ m\leq n \}}
\ee

\noindent
One important consequence of the fact that a process enjoys the self-duality property is that the dynamics of $m$ particles provides relevant information about the time-dependent correlation functions of degree $m$. As an example we now state  the following proposition, Proposition 5.1  in \cite{carinci2017exact}, which provides evidence for the case of two particles

\bp\label{timedepcorr}
Let $\{\eta(t): t\geq 0\}$ be a process with generator \eqref{SIPgen}, then
\beq\label{simplifVariance2}
&& \int \E_{\eta}\( \eta_t(x) - \rho \) \( \eta_t(y) - \rho \) \nu (d \eta)  \\
&&= \( 1 + \frac{1}{k} \1_{ \{ x= y \}} \) \( \frac{ k \sigma}{k +1} - \rho^2 \) \E_{x,y} \left[ \1_{ \{ X_t = Y_t \} } \right]  + \1_{ \{ x= y \}} \( \frac{\rho^2}{k} + \rho \)\nn
\eeq
where $\nu$ is assumed to be a homogeneous product measure with  $\rho$ and $\sigma$  given by
\beq\label{rhoandsigma}
\rho := \int \eta_x \nu ( d \eta ) \qquad \text{and} \qquad
\sigma := \int \eta_x (\eta_x -1) \nu ( d \eta )
\eeq
and $X_t$ and $Y_t$ denote the positions at time $t >0$ of two dual particles started at $x$ and $y$ respectively and $\E_{x,y}$ the corresponding expectation. 
\ep

\bpr
We refer to \cite{carinci2017exact} for the proof.
\epr

\br\label{RelevanceofDiff}
Notice that Proposition \ref{timedepcorr} shows that the two-point correlation functions depend on the  two particles dynamics via the indicator function $\1_{ \{ X_t = Y_t \}}$. More precisely, these correlations can be expressed in terms of the difference of the positions of two dual particles and the model parameters.
\er

\noindent
Motivated by Remark \ref{RelevanceofDiff}, and for reasons that will become clear later, we will study in the next section the stochastic process obtained from the generator \eqref{SIPgen} by following the evolution in time of the difference of the positions of two dual particles.

\subsection{The difference process}

We are interested in a  process obtained from  the dynamics of the process $\{\eta(t):t\geq 0\}$ with generator \eqref{SIPgen} initialized originally with two labeled particles. More precisely, if we denoted by $(x_1(t),x_2(t))$ the particle positions at time $t \geq 0$, from the generator \eqref{SIPgen} we can deduce the generator for the evolution of these two particles; this is, for $f:\Z^2 \to \R$ and $ \mathbf{x} \in \Z^2$ we have
\be
L f(\mathbf{x}) = \sum_{i=1}^{2} \sum_{r} p(r) \bigg( k + \sum_{j=1}^{2} \1_{x^i +r = x^j} \bigg) \( f(\mathbf{x}^{i,r}) - f(\mathbf{x})  \) \nn
\ee
where $\mathbf{x}^{i,r}$ results from changing the position of particle $i$ from the site $x^i$ to the site $x^i +r$.\\

\noindent
Given this dynamics, we are interested in the process given by the difference
\be\label{distdef}
w(t) :=  x_2(t) - x_1(t), \qquad t\ge 0
\ee
Notice that once fixed the initial position of particles, the particles keep the same label. This process was studied for the first time in \cite{opoku2015coupling} and later on \cite{carinci2017exact}, but in contrast to \cite{carinci2017exact}, we do not restrict ourselves to the nearest neighbor case, hence at any Poisson clock ring the value of $w(t)$ can change by $r$ units, with $r \in A:= [-R,R] \cap \Z \setminus \{0\}$.\\

\noindent
Using the symmetry and translation invariance properties of the transition function we obtain the following operator as generator for the difference process
\be\label{wgen}
    (Lf) (w) = \sum_{r \in A } 2 p(r) \left( k + \1_{r=-w} \right) \left[ f(w+r)-f(w)\right]
\ee
where we used that $p(0)=0$ and $p(-r)=p(r)$.
\vskip.2cm
\noindent
Let  $\mu$  denote the discrete counting measure and $\delta_0$ the Dirac measure at the origin, we have the following result
\bp\label{reversinurho}
The difference process is reversible with respect to the measure $\nu_k$ given by
\be\label{nurho}
    \nu_k:=  \mu + \frac{\delta_0}{k}, \qquad \text{i.e.}\qquad \nu_k(w)=
    \begin{cases*}
      1+\frac{1}{k} & \text{if} \: $w = 0$ \\
        1  & \text{if}\:   $w \neq 0$
    \end{cases*}
\ee
\ep

\bpr
By detailed balance, see for example Proposition 4.3 in \cite{kipnis2013scaling},  we  obtain that any reversible measure should satisfy the following:
\be\label{detailnuk}
\nu_k(w) = \dfrac{ \left( k + \1_{w=0} \right)}{ \left( k + \1_{r=-w} \right)} \;\nu_k (w+r)
\ee
where, due to the symmetry of the transition function, we have cancelled the factor $\tfrac{p(-r)}{p(r)}$. In order to verify that $\nu_k$ satisfies \eqref{detailnuk} we have to consider 3 possible cases: $w\neq 0,-r$, $w=0$ and $w=-r$. For $w\neq 0,-r$, \eqref{detailnuk} reads $\nu_k(w)=\nu_k(w+r)$ that is clearily satisfied by \eqref{nurho}. For $w=0$ and for $w=-r$, \eqref{detailnuk} reads $\nu_k(0)=(1+\frac 1k)\nu_k(r)$ that is also satisfied by \eqref{nurho}.
\epr

\br\label{independrangemeasu}
Notice that in the case of a symmetric transition function the  reversible measures $\nu_k$ are independent of the range of the transition function.
\er

\subsection{Condensation and Coarsening}

\subsubsection{The sticky regime}\label{sect:cond}

It has been shown in \cite{grosskinsky2011condensation} that the inclusion process with generator \eqref{SIPgen} can exhibit a condensation transition in the limit of a vanishing diffusion parameter $k$. The parameter $k$ controls the rate at which particles perform random walks, hence in the limit $k \to 0$ the interaction due to inclusion becomes dominant which leads to condensation. The type of condensation in the SIP is different from other particle systems such as zero range processes, see \cite{Grosskinsky2003CondensationIT} and \cite{evans2005nonequilibrium} for example, because in the SIP the critical density is zero. \\ 

\noindent
In the symmetric inclusion process we can achieve the condensation regime  by rescaling  the parameter $k$, i.e. making it of order $1/N$.
If on top of that rescaling we  also rescale space by  $1/N$ and accelerate time with a factor of order $N^3$ then we enter the sticky regime introduced in \cite{carinci2017exact}.  More precisely, for $\gamma >0$,  we speed up time by a factor ${N^3 \gamma}/{\sqrt{2}}$, scale space by $1/N$ and rescale  the parameter $k$ by  $\tfrac{1}{\sqrt{2} \gamma  N}$; in this case the generator \eqref{SIPgen} becomes
\be\label{SIPgencondensive}
\loc_N f(\eta) = \frac{N^3 \gamma}{\sqrt{2}} \sum_{i,j \in \tfrac{1}{N}\Z} p(j-i) \eta_i (\tfrac{1}{\sqrt 2 \gamma  N} +\eta_j) (f(\eta^{ij})-f(\eta))
\ee
Notice that by splitting the generator \eqref{SIPgencondensive} as follows:
\be\label{SIPgencondensivesplit}
\loc_N f(\eta) = \loc_N^{\text{IRW}} f(\eta) + \loc_N^{\text{SIP}} f(\eta) \nn
\ee
where
\be\label{SIPgencondensiveIRW}
\loc_N^{\text{IRW}}  f(\eta) = \frac{N^2 }{2} \sum_{i,j \in \tfrac{1}{N}\Z} p(j-i) \eta_i   (f(\eta^{ij})-f(\eta))
\ee
and
\be\label{SIPgencondensiveSIP}
\loc_N^{\text{SIP}} f(\eta) = \frac{N^3 \gamma}{\sqrt{2}} \sum_{i,j \in \tfrac{1}{N}\Z} p(j-i) \eta_i \eta_j  (f(\eta^{ij})-f(\eta))
\ee
We can indeed see two forces competing with each other. On the one hand, with a multiplicative factor of $\frac{N^2 }{2}$  we see the diffusive action of the generator \eqref{SIPgencondensiveIRW}. While on the other hand, at a much larger factor  $\frac{N^3 \gamma}{\sqrt{2}} $ we see the action of the infinitesimal operator \eqref{SIPgencondensiveSIP} making particles condense. Therefore the sum of the two generators have the flavor of a slow-fast system. This gives us the hint that for the associated process we cannot expect convergence of the generators. Instead, as it will become clear later, we will work with Dirichlet forms.

\subsubsection{Coarsening and the density fluctuation field}

It  was found in \cite{grosskinsky2011condensation} that in the condensation regime ( when started from a homogeneous product measure with density $\rho >0$) sites are either empty with very high probability, or contain a large number of particles to match the fixed expected value of the density. We also know that in this regime the variance of the particle number is of order $N$ and hence a rigorous hydrodynamical description of the coarsening process, by means of standard techniques, becomes inaccessible. Nevertheless, as it was already hinted in \cite{carinci2017exact} at the level of the Fourier-Laplace transform, a rigorous description at the level of fluctuations might be possible. Therefore we introduce the fluctuation field in the the condensive time scaling:
\be\label{flucfield}
\caX_N(\eta,\phi, t)= \frac{1}{N} \sum_{x \in \Z} \phi(x/N)  \( \eta_{\alpha (N, t)}(x) - \rho \)  \qquad \text{with} \qquad \alpha (N, t):= \tfrac{\gamma N^3  t}{ \sqrt{2}}
\ee
defined for any Schwartz function $\phi: \R \to \R$.

\br
Notice that the scaling in \eqref{flucfield} differs from the standard setting of fluctuation fields, given for example in Chapter 11 of \cite{kipnis2013scaling}. In our setting, due to the exploding variances it is necessary to re-scale the fields by an additional factor of $\tfrac{1}{\sqrt{N}}$.
\er

\section{Main result: time dependent variances of the density field}

Let us initialize the nearest neighbor SIP configuration process from a spatially homogeneous product measure $\nu$ parametrized by its mean $\rho$ and such that
\be 
\E_{\nu} [\eta(x)^2 ]< \infty \nn
\ee  
We have the following result concerning the time dependent variances of the density field \eqref{flucfield}:

\bt\label{MainAppThm}
Let $\{ \eta_{\alpha(N,t)} : t \geq 0\}$ be the condensively rescaled inclusion process in configuration space. Consider the fluctuation field $\caX_N(\eta,\phi, t)$ given by \eqref{flucfield}. Let $\nu$ be an initial homogeneous product measure then the time variances of the field are such that
\beq\label{timevarthmRM}
 \lim_{N \to \infty}\E_\nu \left[ \caX_N(\eta,\phi, t)^2 \right] &=& -{\sqrt 2 \gamma^2 \rho^2}\; e^{4 \gamma^2 t} \int_{\R^2}  \phi(x) \phi(y) \; e^{2 \sqrt 2 \gamma \lvert x-y \rvert}  \erf( 2\gamma \sqrt{t} +\tfrac{\lvert x-y \rvert}{\sqrt{2t}}) \, dx \, dy   \nn \\
&&\hskip-1cm+ {\sqrt 2 \gamma \rho^2}\( 1-e^{4 \gamma^2 t} \erf( 2\gamma \sqrt{t}) \)\int_{\R} \phi(x)^2 \, dx
\eeq
where we have used the convention $\erf(x) := \frac{2}{\sqrt{\pi}} \int_x^{\infty} e^{-y^2} dy$.
\et

\subsection{Heuristics of the coarsening process}
In this section we give some intuition about the limiting behavior of the density field, as found in Theorem \ref{MainAppThm}.
More concretely,  we show that Theorem \ref{MainAppThm} is consistent with the following ``coarsening picture''.
From the initial homogeneous product measure $\nu$ with density $\rho$, in the course of time large piles are created which are typically at distances of order $N$ and of size $\rho N$.
The location of these piles evolves on the appropriate time scale according to a diffusion process. If we focus on two piles, this diffusion process is of the form $(X(t), Y(t))$ where
$X(t)-Y(t)$ is a sticky Brownian motion $B^{\text{sbm}}(t)$, and where the sum $X(t)+Y(t)$ is an independent Brownian motion $\overline{B}(t)$, time changed via the local time inverse at the origin $\tau(t)$ of the sticky Brownian motion $B^{\text{sbm}}(t)$ via
$X(t)+ Y(t)= \overline{B} (2t- \tau(t))$.\\

\noindent
Let us now make this heuristics more precise. Define the non-centered field
\be\label{flucfield1}
\caZ_N(\eta,\phi, t)= \frac{1}{N} \sum_{x \in \Z} \phi(\tfrac xN)  \eta_{\alpha (N, t)}(x)
\ee
then one has, using that at every time $t>0$, and $x\in\Zd$, $\E_\nu (\eta_t(x))=\rho$:
\be\label{bam}
 \lim_{N \to \infty}\E_\nu \left[ \caZ_N(\eta,\phi, t) \right] = \rho\int_{\R}  \varphi(x)\, dx
\ee
and
\be
 \lim_{N \to \infty}\left(\E_\nu \left[ \caZ_N(\eta,\phi, t)^2 \right]- \E_\nu \left[ \caX_N(\eta,\phi, t)^2 \right]\right) = \rho^2\int_{\R}  \int_{\R}   \varphi(x)\varphi(y)\, dx\, dy\nn
\ee
as we will see later in the proof of our main theorem, the RHS of \eqref{timevarthmRM} can be written as \eqref{quo} and hence we have that
\beq\label{bim}
&& \lim_{N \to \infty}\E_\nu \left[ \caZ_N(\eta,\phi, t)^2 \right] =\rho^2\int_{\R} \int_{\R}  \varphi(x)\varphi(y)\, dx\, dy+ \sqrt 2 \gamma \rho^2  \int_{\R} \phi(u)^2 \, du\nn \\
&& -\frac{\rho^2}{2} \int_{\R^2}  \phi(\tfrac{u+v}{2}) \phi(\tfrac{u-v}{2}) p_t^{\text{sbm}}(v,0) \, dv \, du   - \sqrt 2 \gamma \rho^2 p_t^{\text{sbm}}(0,0) \int_{\R} \phi(u)^2 \, du  \nn \\
&& =\frac{\rho^2}2\int_{\R} \int_{\R}  \varphi\(\tfrac{u+v}2\)\varphi\(\tfrac{u-v}2\)\,(1-p_t^{\text{sbm}}(v,0))\, (dv+ \sqrt 2 \gamma\delta_0(dv))\, du\nn\\
&&= \frac{\rho^2}{2} \int_{\R} \( \int_{\R} \E_v^{\text{sbm}}  \( \phi(\tfrac{u+v_t}{2}) \phi(\tfrac{u-v_t}{2}) \) (1-\1_{\{ 0 \}}(v)) \( \, dv   + \sqrt 2 \gamma \delta_0(dv)\) \) \, du \nn \\
&&= \frac{\rho^2}{2} \int_{\R} \int_{\R} \E_v^{\text{sbm}}  \( \phi(\tfrac{u+v_t}{2}) \phi(\tfrac{u-v_t}{2}) \)   \, dv\, du \nn\\
&&= \frac{\rho^2}{2} \int_{\R}\int_{\R} \int_{\R} \phi(\tfrac{u+z}{2}) \phi(\tfrac{u-z}{2}) \,  p_t^{\text{sbm}}(v,dz)   \, dv\, du \nn\\
&&= {\rho^2}\int_{\R}dv \int_{\R}\int_{\R}  \phi(x) \phi(y) \cdot \bar p_t^{\text{sbm}}(v;dx,dy)
\eeq
where
\be
\bar p_t^{\text{sbm}}(v;dx,dy):= p_t^{\text{sbm}}(v,x-y) \, dx\, dy+p_t^{\text{sbm}}(v,0) \, dx\, \delta_x(dy)
\ee
 in the second line we used the change of variables
$x = \tfrac{u+v}{2}$, $y = \tfrac{u-v}{2}$. 
\noindent
We now want to describe a ``macroscopic'' time dependent random field $\caZ(\phi, t)$ that is consistent with the limiting expectation and second moment
computed in \eqref{bam} and \eqref{bim}. This macroscopic field describes intuitively the positions of the piles formed from the initial homogeneous background.\\

\noindent
First define for $m\in \N$
\be
\caZ^{(m)}(\phi, t)= \frac{\rho}m\sum_{i=1}^m\int_{\R} \phi(X_i^x(t)) dx
\ee
where $(X^{x_1}_1(t),\ldots,X^{x_m}_m(t))$ is a $m$-dimensional diffusion process such that
\begin{itemize}
\item[a)] the marginals $X^{x_i}_i(t)$ are  Brownian motions with diffusion constant $\chi/2 $ started from $x_i$.
\item[b)] the couples $\{(X^{x_i}_i(t),X^{x_j}_j(t)), \; i,j=1,\ldots m\}$  are two dimensional diffusion processes starting from initial positions $(x_i,x_j)$. At any fixed time $t\ge 0$ each couple is distributed in such a way that the difference-sum process  is given by
\be\label{couple}
(X^{x_i}_i(t)-X^{x_j}_j(t)),X^{x_i}_i(t)+X^{x_j}_j(t))= (B^{\text{sbm}, x_i-x_j}(t),\bar B^{x_i+x_j}(2t-\tau(t)))
\ee
Here $B^{\text{sbm}, x_i-x_j}(t)$ is a sticky Brownian motion with stickiness parameter $\sqrt 2 \gamma$, and diffusion constant $\chi$, started from $x_i-x_j$ and where $\tau(t)$ is the corresponding local time-change defined in \eqref{tauinv}, and $\bar B^{x_i+x_j}(2t-\tau(t))$ is another Brownian motion and diffusion constant $\chi$, independent from $B^{\text{sbm}}(t)$ started from
$x_i+x_j$.

\end{itemize}
Then we will see that in the limit $m\to\infty$, the field $\caZ^{(m)}(\phi, t)$ reproduces correctly the first and second moments
of \eqref{bam} and \eqref{bim}.\\

\noindent
For the expectation we have, using item a) above
\beq
\mathbb E[\caZ^{(m)} (\phi, t)]=\frac{\rho}m\sum_{i=1}^m  \int_{\R} \mathbb E[\phi(X^x_i(t))] dx_i=\rho \int_{\R} \phi(x)\int_{\R}  p_t^{\text{bm}}(x_0,x)\, dx_0\, dx=\rho \int_{\R}  \phi(x)\, dx\nn
\eeq
where the last identity follows from the symmetry: $p_t^{\text{bm}}(x_0,x)=p_t^{\text{bm}}(x,x_0)$.\\

\noindent
On the other hand, for the second moment, using item b) above
\beq
\mathbb E[\caZ^{(m)}(\phi, t)^2]= \frac{\rho^2}{m^2}\sum_{i,j=1}^m\int_{\R}\int_{\R}  \mathbb E[\phi(X^{x_0}_i(t))\phi(X^{y_0}_j(t))] dx_0 dy_0
\eeq
Let $i\neq j$, then, from our assumptions,
\beq
&&\mathbb E[\phi(X^{x_0}_i(t))\phi(X_j^{y_0}(t))] =   \int_{\R}  \int_{\R} \phi(x)\phi(y) p_t(x_0,y_0;dx,dy) \nn
\eeq
Where  $p_t(x_0,y_0;dx,dy)$ is the transition probability kernel  of the couple $(X_1(t),X_2(t))$. Denoting now by $\tilde p_t(v_0,u_0;dv,du)$ the transition probability kernel of the couple $(X_1(t)-X_2(t),X_1(t)+X_2(t))$,  and by $\pi_t$ the probability measure of the time change $\tau(t)$, at time $t$. Then we have
\beq
 \tilde p_t(v_0,u_0;dv,du)= \int_{\R}   \tilde p_t(v_0,u_0;dv,du \, | s)\, \pi_t(ds) = \int_{\R}    \tilde p^{(1)}_t(v_0,dv \, | s)\,  \tilde p^{(2)}_t(u_0,du \, | s)\,\pi_t(ds) \nn
\eeq
(where $ \tilde p^{(i)}_t(\cdot,\cdot| s)$ for $i=1,2$, are resp. the transition probability density functions of the Brownian motions $B(t)$ and $\bar B(t)$ conditioned on $s$)  as, from \eqref{couple}, the difference and sum processes are independent conditionally on the realization of $s(t)$.
Now we have that
\beq
\int_{\R}  \tilde p^{(1)}_t(v_0,dv \, | s) \pi_t(ds)= p^{\text{sbm}}_{t}(v_0,dv)
 \quad \text{and} \quad \tilde p^{(2)}_t(u_0,du \, | s)= p^{\text{bm}}_{2t-s}(u_0,du) \nn
\eeq
hence
\beq
\int_{\R}  \int_{\R}  \tilde p_t(v_0,u_0;dv,du) \,dv_0\,du_0&= &\int_{\R}    \(\int_{\R}  \tilde p^{(1)}_t(v_0,dv \, | s)\,dv_0 \) \cdot  \(\int_{\R}    p^{\text{bm}}_{2t-s}(u_0,du)\, du_0 \)\,\pi_t(ds) \nn\\
& =& \int_{\R}    \int_{\R}  \tilde p^{(1)}_t(v_0,dv \, | s) \pi_t(ds) \, dv_0 = \int_{\R}      p^{\text{sbm}}_t(v_0,dv) \, dv_0
\eeq
where the second identity follows from the symmetry of $p^{\text{bm}}(\cdot,\cdot)$. Then,
 from the change of variables $v_0:=x_0-y_0$, $u_0=x_0+y_0$, and $v=x-y$, $u=x+y$, and since $dv_0\,du_0=2 dx_0\,dy_0$, 
it follows that
\beq\int_{\R}  \int_{\R}  p_t(x_0,y_0;dx,dy)\,dx_0\,dy_0=  \int_{\R}     \bar p_t^{\text{sbm}}(v_0;dx,dy) \, dv_0
\eeq
For $i=j$ we have
\beq
&&\mathbb E[(\phi(X_i^{x_0}(t)))^2]= \int_{\R}   (\phi(x))^2    p^{\text{bm}}_t(x_0,dx)   \nn
\eeq
then
\beq
&&\mathbb E[(\caZ^{(m)}(\phi, t))^2]\nn\\
&&=\rho^2 \int_{\R}  \int_{\R}  \phi(x)\phi(y) \int_{\R}  \left\{ \(1-\frac 1 m\) \,    \bar p_t^{\text{sbm}}(v;dx,dy)  + \frac 1 m  \, p^{\text{bm}}_t(v,dx)  \delta_{x}(dy)\right\} \, dv \nn
\eeq
this converges to
\beq
\rho^2 \int_{\R}  \int_{\R}  \phi(x)\phi(y) \int_{\R}   \bar p_t^{\text{sbm}}(v;dx,dy) \, dv 
\eeq
in the limit as $m \to \infty$.

\section{Basic tools}

Before showing the main result, in this section we introduce some notions and tools that will be useful to show Theorem \ref{MainAppThm}. These notions include the concept of Dirichlet forms and the notion of convergence of Dirichlet forms that we will use; Mosco convergence of Dirichlet forms. The reader familiar with these notions can  skip this section and move directly to Section \ref{proofs}.

\subsection{Dirichlet forms}

A Dirichlet form on a Hilbert space is defined as a symmetric form which is closed and Markovian. The importance of Dirichlet forms in the theory of Markov processes is that the Markovian nature of the first corresponds to the Markovian properties of the associated semigroups and resolvents on the same space. Related to the present work, probably one of the best examples of this connection is the work of Umberto Mosco. In \cite{mosco1994composite} Mosco introduced a type of convergence of quadratic forms, Mosco convergence, which is equivalent to strong convergence of the corresponding semigroups. Before defining this notion of convergence, we recall the precise definition of Dirichlet form.

\bd[Dirichlet forms]
Let $H$ be a Hilbert space of the form $L^2(E;m)$ for some $\sigma$-finite measure space $(E,\caB(E),m)$. Let $H$ be endowed with an inner product $ \langle \cdot,\cdot \rangle_H$. A Dirichlet form $\caE (f,g)$ on $H$ is a symmetric bilinear form such that the following conditions hold
\ben
\item The form is closed, i.e. the domain $ D(\caE)$ is complete with respect to the metric determined by
\be
\caE_1 (f,g) = \caE (f,g) +  \langle f,g \rangle_H \nn
\ee
\item The unit contraction operates on $\caE$, i.e. for $ f \in D(\caE)$, if we set $ g := ( 0 \vee f ) \wedge 1$ then we have that $g \in D(\caE)$ and $\caE(g,g) \leq \caE(f,f)$.
\een
\ed
\noindent
When the second condition is satisfied we say that the form $\caE$ is Markovian. We refer the reader to  \cite{fukushima1980dirichlet} for a comprehensible  introduction to the subject of Dirichlet forms. For the purposes of this work, the key property of Dirichlet forms is that there exists a natural correspondence between the set of Dirichlet forms and the set of Markov generators. In other words, to a symmetric Markov process we can always associate a Dirichlet form that is given by:
\be\label{DirGenRel}
\caE(f,g)= - \langle f,Lg \rangle_H \qquad \quad \text{with} \qquad
D(\caE)= D(\sqrt{-L})
\ee
where the operator $L$ is the corresponding infinitesimal generator of a symmetric Markov process. As an example of this relation, consider the  Brownian motion in $\R$. We know that the associated infinitesimal generator is given by the Laplacian. Hence its Dirichlet form is
\beq\label{SBMintro}
\caE_{\text{bm}}(f, g) = \frac{1}{2} \int_{-\infty}^\infty f\myprime(x) g\myprime (x) dx
\qquad \text{with domain} \qquad
D(\caE_{\text{bm}})= H^1(\R)
\eeq
namely the Sobolev space of order 1.\\

\noindent
From now  on  we will mostly deal with the quadratic form $ \caE(f,f)$ that we can view as a functional defined on  the entire Hilbert space $H$ by defining
\be\label{ext2H}
    \caE(f) =
    \begin{cases*}
      \caE(f,f), &  $f \in D(\caE) $ \\
      \infty,  & $f \notin D(\caE), $
    \end{cases*}
\qquad f\in H
\ee
which is lower-semicontious if and only if the form $(\caE, D(\caE)) $ is closed. \\

\subsection{Mosco convergence}

We now introduce the framework to properly define the mode of convergence we are interested in. The idea is that we want to approximate a Dirichlet form on the continuum by a sequence of Dirichlet forms indexed by a scaling parameter $N$. In this context, the problem with the convergence introduced in \cite{mosco1994composite} is that the approximating sequence of Dirichlet forms does not necessarily live on the same Hilbert space. However, the work in \cite{kuwae2003convergence} deals with this issue. We also refer to \cite{kolesnikov2006mosco} for a more complete understanding and a further generalization to infinite dimensional spaces. 
In order to introduce this mode of convergence, we  first define the concept of convergence of Hilbert spaces.

\subsection{Convergence of Hilbert spaces}
We start with the definition of the notion of convergence of spaces:

\bd[Convergence of Hilbert spaces]\label{HilConv}
A sequence of Hilbert spaces $\{ H_N \}_{N \geq 0}$, converges to a Hilbert space $H$ if there exist a dense subset $C \subseteq H$ and a family of linear maps $\{ \Phi_N : C  \to H_N \}_N$ such that:
\be\label{HilCond}
\lim_{N \to \infty} \| \Phi_N f \|_{H_N} = \| f \|_{H}, \qquad \text{  for all } f \in C
\ee
\ed
\noindent
It is also necessary to introduce the concepts of strong and weak convergence of vectors living on a convergent sequence of Hilbert spaces. Hence in Definitions \ref{strongcon}, \ref{weakcon} and \ref{MoscoDef} we assume that the spaces $\{ H_N \}_{N \geq 0}$ converge to the space $H$, in the sense we just defined,  with the dense set  $C \subset H$ and the sequence of operators $\{ \Phi_N : C  \to H_N \}_N$ witnessing the convergence.

\bd[Strong convergence on Hilbert spaces]\label{strongcon}
A sequence of vectors $\{  f_N \}$ with $f_N$ in $H_N$, is said to strongly-converge to a vector $f \in H$, if there exists a sequence $\{ \tilde{f}_M \} \in C$ such that:
\be
\lim_{M\to \infty} \| \tilde{f}_M -f \|_{H} = 0
\ee
and
\be
\lim_{M \to \infty} \limsup_{N \to \infty} \| \Phi_N \tilde{f}_M -f_N \|_{H_N} = 0
\ee
\ed
\bd[Weak convergence on Hilbert spaces]\label{weakcon}
A sequence of vectors $\{  f_N \}$ with $f_N \in H_N$, is said to converge weakly to a vector $f$ in a  Hilbert space $H$ if
\be
\lim_{N\to \infty} \left \langle f_N, g_N \right \rangle_{H_N} = \ \left \langle f, g \right \rangle_{H}
\ee
for every sequence $\{g_N \}$ strongly convergent to $g \in H$.
\ed
\br
Notice that, as expected, strong convergence implies weak convergence, and, for any $f \in C$, the sequence $\Phi_N f $ strongly-converges to $f$.
\er

\noindent
Given these notions of convergence, we can also introduce related notions of convergence for operators. More precisely, if we denote by $L(H)$ the set of all bounded linear operators in $H$, we  have the following definition

\bd[Convergence of bounded operators on Hilbert spaces]\label{opercon}
A sequence of bounded operators $\{  T_N \}$ with $T_N \in L(H_N)$, is said to strongly (resp. weakly ) converge to an operator  $T$ in $L(H)$ if for every strongly (resp. weakly) convergent sequence $\{  f_N \}$, $f_N \in H_N$ to $f \in H$ we have that the sequence $\{ T_N f_N \}$ strongly (resp. weakly ) converges to $T f$.
\ed

\noindent
We are now ready to introduce Mosco convergence.

\subsection{Definition of Mosco convergence}
In this section we assume the Hilbert convergence of a sequence of Hilbert spaces $\{ H_N \}_N$ to a space $H$.

\bd[Mosco convergence]\label{MoscoDef}
A sequence of Dirichlet forms $\{ (\caE_N, D(\caE_N))\}_N $ on Hilbert spaces $H_N$, Mosco-converges to a Dirichlet form $(\caE, D(\caE)) $ in some Hilbert space $H$ if:
\begin{description}
\item[Mosco I.] For every sequence of $f_N \in H_N$ weakly-converging  to $f$ in $H$
\be\label{mosco1}
\caE ( f ) \leq \liminf_{N \to \infty} \caE_N ( f_N )
\ee
\item[Mosco II.] For every $f \in H$, there exists a sequence $ f_N \in H_N$ strongly-converging  to $f$ in $H$, such that
\be
\caE ( f) = \lim_{N \to \infty} \caE_N ( f_N )
\ee
\end{description}
\ed

\noindent
As we mentioned before, the Markovian properties of the Dirichlet form correspond to the properties of the associated semigroups and resolvents. The following theorem from \cite{kuwae2003convergence}, which relates Mosco convergence with convergence of semigroups and resolvents, is a powerful application of this correspondence and one of the main ingredients of our work:
\bt\label{MKS}
Let $\{ (\caE_N, D(\caE_N))\}_N $ be a sequence of Dirichlet forms on Hilbert spaces $H_N$ and let  $(\caE, D(\caE)) $ be a Dirichlet form  in some Hilbert space $H$. The following statements are equivalent:
\ben
\item $\{ (\caE_N, D(\caE_N))\}_N $ Mosco-converges to $\{ (\caE, D(\caE))\} $.
\item The associated sequence of semigroups $\{ T_{N} (t) \}_N $ strongly-converges to the semigroup $ T(t)$ for every $t >0$.
\een
\et

\subsection{Mosco convergence and dual forms}\label{simpli}
The difficulty in proving condition Mosco I lies in the fact that \eqref{mosco1}  has to hold for all weakly convergent sequences, i.e., we cannot choose a particular class of sequences. \\

\noindent
In this section we will show how  one can avoid this difficulty by passing to the dual form. We prove indeed   that  Mosco I for the original form is implied by a condition similar to Mosco II for the dual form (Assumption \ref{strongcontinuity}).

\subsubsection{Mosco I}\label{simpliMI}

Consider a sequence of Dirichlet  forms $( \caE_N, D(\caE_N) )_N$ on Hilbert spaces $H_N$, and an additional quadratic form $(\caE, D(\caE))$ on a Hilbert space $H$. We assume
 convergence of Hilbert spaces, i.e. that there exists a dense set $C\subset H$  and a sequence of maps
$\Phi_N: C\to H_N$ such that $\lim_{N\to \infty} \|\Phi_N f\|_{H_N} =\|f\|_H$.
The dual quadratic form is defined via
\[
\caE^* (f)= \sup_{g \in H} \( \left \langle f,g \right \rangle- \caE(g)\)
\]
Notice that from the convexity of the form we can conclude that it is involutive, i.e.,  $(\caE^*)^*=\caE$. We now assume that the following holds
\basu\label{strongcontinuity}
For all $ g \in H$, there exists a sequence $ g_N \in H_N $  strongly-converging to $g$ such that
\be\label{dualcondual}
\lim_{N\to \infty}\caE^*_N(g_N)= \caE^*(g)
\ee
\easu
\noindent
We show now that, under Assumption \ref{strongcontinuity}, the first condition of Mosco convergence is satisfied.

\bp
 Assumption \ref{strongcontinuity} implies Mosco I, i.e.
\be
 \liminf_{N \to \infty} \caE_N ( f_N ) \ge \caE ( f )
\ee
for all $f_N \in H_N $ weakly-converging  to $f \in H$.
\ep
\bpr
Let $f_N\to f$ weakly then, by Assumption \ref{strongcontinuity}, for any $g\in H$ there exists a sequence $g_N \in H_N$ such that $  g_N\to g$ strongly, and \eqref{dualcondual} is satisfied. Fromt the involutive nature of the form, and by Fenchel's inequality, we obtain:
\be
\caE_N(f_N) = \sup_{h\in H_N} \(\left \langle f_N, h \right \rangle_{H_N}- \caE^*_N(h)\) \geq \left \langle f_N, g_N \right \rangle_{H_N}- \caE^*_N(g_N) \nn
\ee
by the fact that $f_N \to f$ weakly, $  g_N\to g$ strongly, and \eqref{dualcondual} we obtain
\begin{eqnarray*}
\liminf_{N \to \infty} \caE_N(f_N)\geq \liminf_{N \to \infty} \(\left \langle f_N, g_N \right \rangle_{H_N}- \caE^*_N(g_N)\)
\geq
\left \langle f,g \right \rangle_H -\caE^* (g)
\end{eqnarray*}
Since this holds for all $g \in H$ we can take the supremum over $H$,
\beq\label{supremumdual}
\liminf_{N \to \infty} \caE_N (f_N) \geq \sup_{g \in H} \( \left \langle f,g \right \rangle_H -\caE^* (g)\)=\caE(f)
\eeq
This concludes the proof.
\epr

\noindent
In other words, in order to prove condition Mosco I all we have to show is that Assumption \ref{strongcontinuity} is satisfied.

\subsubsection{Mosco II}

For the second condition we recall a result from \cite{andres2010particle} in which
 a weaker notion of Mosco convergence is proposed, where Mosco I is unchanged whereas Mosco II is relaxed
 to functions living in a core of the limiting Dirichlet form:
\basu\label{secondass}
There exists a core $K \subset D(\caE)$  of $\caE$ such that, for every $f \in K$, there exists a sequence $\{ f_N \}$  strongly-converging to $f$, such that
\be\label{equasecondass}
\caE ( f) = \lim_{N \to \infty} \caE_N ( f_N )
\ee
\easu
\noindent
Despite of being weaker, the authors were able to prove that this relaxed notion also implies strong convergence of semi-groups.  We refer the reader to Section 3 of \cite{andres2010particle} for details on the proof.

\section{Proof of main result}\label{proofs}

Our main theorem, Theorem \ref{MainAppThm}, is a consequence of self-duality and the following result concerning the convergence in the Mosco sense of the sequence of Dirichlet forms associated to the difference process to the Dirichlet form corresponding to the so-called two sided sticky Brownian motion (See the Appendix for details on this process):

\bt\label{MainThm}
The sequence of Dirichlet forms $\{ \caE_N, D(\caE_N)\}_{N \geq 1}$ given by \eqref{DirwN} converges in the Mosco sense to the form $\( \caE_{\text{sbm}}, D(\caE_{\text{sbm}}) \)$ given by
\be\label{DIRSBMchi}
\caE_{\text{sbm}} (f) = \frac \chi 2  \int_{\R} f\myprime(x)^2 dx, \qquad \chi=\sum_{r=1}^Rr^2\, p(r)
\ee
whose  domain is
\be\label{DomSBMchi}
 D(\caE_{\text{sbm}}) =H^1(\R) \cap L^2(\R, \bar\nu) \qquad \text{with} \qquad \bar\nu=dx + \sqrt 2 \gamma \delta_0
\ee
As a consequence, if we denote by $T_N(t)$ and $T_t$  the semigroups associated to the difference process $w_N(t)$ and the sticky Brownian motion $B_t^{\text{s}}$, we have that $T_N(t) \to T_t$ strongly in the sense of Definition \ref{opercon}.
\et
We will show in the following section, how to use this result to prove Theorem Theorem \ref{MainAppThm}. The proof of Theorem \ref{MainThm} will be left to Section \ref{MCInclusec}.

\subsection{Proof of main theorem: Theorem \ref{MainAppThm}}
We then denote by $T_N(t)$ and $T_t$  the semigroups associated to the difference process $w_N(t)$ and the sticky Brownian motion $B_t^{\text{s}}$. Because of our result on Mosco convergence and thanks to Theorem \ref{MKS} we know that the sequence of semigroups $\{ T_N(t) \}_{N \geq 1}$ converges strongly to $T_t$ in the $H_N^{\text{sip}}$ Hilbert convergence sense.  We will see that this implies the convergence of the probability mass function at 0.
\vskip.1cm
\noindent
In the following we denote by $p_t^{\text{sbm}}(x,dy)$ the transition kernel of a Sticky Brownian motion with stickiness parameter $\sqrt 2 \gamma$. This kernel consists of a first term that is absolutely continuous w.r.t. the Lebesgue measure and a second term that is a Dirac-delta at the origin times the probability mass function at zero. With a slight abuse of notation we will denote by 
\be\label{Ciaoo}
p_t^{\text{sbm}}(x,dy)= p_t^{\text{sbm}}(x,y)\, dy + p_t^{\text{sbm}}(x,0) \cdot \delta_0(dy) 
\ee
where $p_t^{\text{sbm}}(x,y)$ for $y\neq 0$ denotes a probability density to arrive at $y$ at time $t$ when started from $x$ , and for $y=0$ the probability to arrive at zero when started at $x$. See equation (2.15) in \cite{howitt2007stochastic} for an explicit formula of \eqref{Ciaoo}. \\

\vskip.1cm
\noindent
We have the following result.
\bp\label{transconv}
For all $t >0$ denote by $p_{t}^N (w, 0)$ the trasition function that the difference process starting from $w \in \frac{1}{N} \Z$ finishes at $0$ at time $t$. Then the sequence $p_{t}^N (\cdot, 0)$ converges strongly to $p_t^{\text{sbm}}(\cdot,0)$ with respect to $H_N^{\text{sip}}$ Hilbert convergence.
\ep
\bpr
From the fact that $\{ T_N(t) \}_{N \geq 1}$ converges strongly to $T_t$, we have that for all $f_N$ strongly converging to $f$, the sequence $\{ T_N(t) f_N \}_{N \geq 1} \in  H_N^{\text{sip}}$ converges strongly to $T_t f$. In particular, for $f_N = \1_{\{0\}}$ we have that the sequence
\be
 T_N(t) f_N (w) = \E_w^N \1_{\{0\}} (w_t) = p_{t}^N (w, 0)
\ee
converges strongly to
\be
T_t f(w) = \E_w^{\text{ sbm} } \1_{\{0\}} (w_t) = p_t^{\text{sbm}}(w,0)
\ee
where $\E_w^{\text{sbm}}$ denotes expectation with respect to the sticky Brownian motion started at $w$.
\epr

\br\label{transconvzero}
Despite the fact that Proposition \ref{transconv} is not a point-wise statement, we can still say something more relevant when we start our process at the point zero:
\be
\lim_{N\to \infty} p_{t}^N (0, 0) = p_t^{\text{sbm}}(0,0)
\ee
The reason is that we can see $p_{t}^N (w, 0)$ as a weakly converging sequence and used again the fact that $f_N = \1_{\{0\}}$ converges strongly.
\er

\bpr{Theorem \ref{MainAppThm}}
Let  $\rho$  and $\sigma$ be given by \eqref{rhoandsigma}, then we can write
\beq
&& \E_\nu \left[ \caX_N(\eta,\phi, t)^2 \right] =  \frac{1}{N^{2}} \sum_{x, y \in \Z} \phi(\tfrac xN) \phi(\tfrac yN) \int \E_{\eta}\( \eta_{\alpha (N, t)}(x) - \rho \) \( \eta_{\alpha (N, t)}(y) - \rho \) \nu (d \eta) \nn 
\eeq
where, from Proposition 5.1 in  \cite{carinci2017exact}, using self-duality we can simplify the integral above as
\beq\label{simplifVariance2}
&& \int \E_{\eta}\( \eta_{\alpha (N, t)}(x) - \rho \) \( \eta_{\alpha (N, t)}(y) - \rho \) \nu (d \eta)  \\
&&= \( 1 + \frac{1}{k_N} \1_{ \{ x= y \}} \) \( \frac{ k_N \sigma}{k_N+1} - \rho^2 \) \E_{x,y} \1_{ \{ X_{\alpha (N, t)} = Y_{\alpha (N, t)} \} }  + \1_{ \{ x= y \}} \( \frac{\rho^2}{k_N} + \rho \)\nn
\eeq
\noindent
Notice that the expectation in the RHS of \eqref{simplifVariance2} can be re-written in terms of our difference process as follows:
\be
\E_{x,y}\left[\1_{ \{ X_{\alpha (N, t)} = Y_{\alpha (N, t)} \} }\right]= p_{\alpha (N, t)} (x-y, 0)
\ee
where $p_{\alpha (N, t)}$ is the transition function $p_t^N$ under the condensive time rescaling defined in \eqref{flucfield}. Since under the condensation regime we have, as in Section \ref{sect:cond}, 
$k_N  = \frac 1 {\sqrt 2 \gamma  N}$. We then obtain:
\beq\label{simplifVariance3}
&& \E_\nu \left[ \caX_N(\eta,\phi, t)^2 \right] \nn \\
&&= \frac{1}{N^{2}} \sum_{x, y \in \Z} \phi(\tfrac xN) \phi(\tfrac yN) \( 1 + \sqrt 2 \gamma  N \1_{ \{ x= y \}} \) \( \frac{ \sigma}{1 +\sqrt 2 \gamma  N} - \rho^2 \) p_{\alpha (N, t)} (x-y, 0)  \nn \\
&&+ \frac{1}{N^{2}} \sum_{x \in \Z} \phi(\tfrac xN) \phi(\tfrac xN) \( \sqrt 2 \gamma  N \rho^2 + \rho \)  
\eeq
\noindent
At this point we have 3 non vanishing contributions:
\beq
&&C_N^{(1)}:= \frac{\rho^2}{N^{2}} \sum_{x, y \in \Z} \phi(\tfrac xN)\phi(\tfrac yN)  p_{\alpha (N, t)} (x-y, 0),  \nn \\
&&\hskip1cm C_ N^{(2)}:=  \frac{\sqrt 2 \gamma\rho^2}{N} \sum_{x \in \Z}\(\phi(\tfrac xN)\)^2  p_{\alpha (N, t)} (0, 0)  \quad \text{and} \qquad C_N^{(3)}:= \frac{\sqrt 2 \gamma \rho^2}{N} \sum_{x \in \Z} \(\phi(\tfrac xN)\)^2 \nn
\eeq
where we already know:
\be\label{C3limit}
\lim_{N\to \infty} C_N^{(3)} =  \sqrt 2 \gamma \rho^2\int_{\R} \phi(v)^2 dv
\ee
and, by Remark \ref{transconvzero},
\be\label{C2limit}
\lim_{N\to \infty} C_ N^{(2)} =  \sqrt 2 \gamma \rho^2 p_t^{\text{sbm}}(0,0) \int_{\R} \phi(v)^2 dv.
\ee
To analyze the first contribution we use the change of variables $u = x+y$, $v = x-y$ from which  we obtain:
\be
C_N^{(1)} = \frac{\rho^2}{N^{2}} \sum_{\substack{u, v \in \tfrac{1}{N}\Z \\ u  \equiv  v \mod 2}} \phi(\tfrac{u+v}{2}) \phi(\tfrac{u-v}{2}) \, p_{\alpha (N, t)} (v, 0)
\ee
hence by \eqref{L2NSip}, $C_N^{(1)}$  can be re-written as
\be
C_N^{(1)} = \left \langle F_N(\cdot), p_{\alpha (N, t)} (\cdot, 0) \right \rangle_{H_{N}^{\text{sip}}} - \frac{ \gamma\rho^2}{\sqrt 2 N} \sum_{u \in \tfrac{1}{N} \Z} \phi(\tfrac{u}{2}) \phi(\tfrac{u}{2}) \, p_{\alpha (N, t)} (0, 0)
\ee
with $F_N$ given by
\be
F_N (v) = \frac{\rho^2}{N} \sum_{\substack{u \in \tfrac{1}{N}\Z \\ u\equiv v\mod 2}} \phi(\tfrac{u+v}{2}) \phi(\tfrac{u-v}{2}), \qquad  \text{ for all } v \in \tfrac{1}{N} \Z
\ee
notice that $F_N$ converges strongly to the function $F \in  H^{\text{sbm}}$ given by
\be
F(x) := \frac{\rho^2}{2} \int_{\R} \phi(\tfrac{y+x}{2}) \phi(\tfrac{y-x}{2}) dy
\ee
which can be seen, in the language of Definition \ref{strongcon}, by setting the reference sequence of functions  $\tilde{F}_M=F$ for all $M$, and from the convergence
\be
\lim_{N \to \infty} \sum_{v \in \tfrac{1}{N} \Z} \int_{\R} \phi(\tfrac{u+v}{2}) \phi(\tfrac{u-v}{2}) du = \int_{\R^2}  \phi(\tfrac{u+v}{2}) \phi(\tfrac{u-v}{2}) \, du \, dv
\ee

\noindent
From the strong convergence $F_N \to F$, Proposition \ref{transconv}, and Remark \ref{transconvzero} we conclude
\be
 \lim_{N \to \infty} C_N^{(1)} = \frac{\rho^2}{2} \int_{\R^2}  \phi(\tfrac{u+v}{2}) \phi(\tfrac{u-v}{2}) p_t^{\text{sbm}}(v,0) \, du \, dv
\ee
substituting the limits of the contributions we obtain
\beq\label{quo}
&& \lim_{N \to \infty}\E_\nu \left[ \caX_N(\eta,\phi, t)^2 \right] \nn \\
&&= -\frac{\rho^2}{2} \int_{\R^2}  \phi(\tfrac{u+v}{2}) \phi(\tfrac{u-v}{2}) p_t^{\text{sbm}}(v,0) \, dv \, du   - \(\sqrt 2 \gamma \rho^2 p_t^{\text{sbm}}(0,0) - \sqrt 2 \gamma \rho^2\)\int_{\R} \phi(u)^2 \, du  \nn \\
&&= -\frac{\rho^2}{2} \int_{\R}  \int_{\R}  \phi(\tfrac{u+v}{2}) \phi(\tfrac{u-v}{2})  \E_v^{\text{sbm}}  \left[\1_{\{ 0 \}}(v_t) \right] \big(dv   + \sqrt 2 \gamma \delta_0(dv)\big)\, du + \sqrt 2 \gamma \rho^2  \int_{\R} \phi(u)^2 \, du \nn \\
&&= -\frac{\rho^2}{2} \int_{\R}  \int_{\R} \E_v^{\text{sbm}}  \left[ \phi(\tfrac{u+v_t}{2}) \phi(\tfrac{u-v_t}{2}) \right] \1_{\{ 0 \}}(v) \big(dv   + \sqrt 2 \gamma \delta_0(dv)\big)\, du + \sqrt 2 \gamma \rho^2  \int_{\R} \phi(u)^2 \, du \nn \\
&&= \frac{\sqrt 2 \gamma \rho^2 }{2}\int_{\R} \left\{ \phi(\tfrac{u}{2})^2 - \E_0^{\text{sbm}}  \left[ \phi(\tfrac{u+v_t}{2}) \phi(\tfrac{u-v_t}{2}) \right] \right\}  \, du\nn\\
&&= \frac{\sqrt 2 \gamma \rho^2 }{2}\int_{\R}\left\{ \phi(\tfrac{u}{2})^2 -   \int_{\R} p_t^{\text{sbm}}(0,dv) \phi(\tfrac{u+v}{2}) \phi(\tfrac{u-v}{2}) \right\}  \, du \label{ciao}
\eeq
where in the third equality we used the reversibility of SBM with respect to the measure $ \hat{\nu}(dv)=dv + \sqrt 2 \gamma \delta_0(dv)$. Then, \eqref{timevarthmRM} follows, after a change of variable, using the expression (2.15) given in \cite{howitt2007stochastic} for the transition probability measure $p_t^{\text{sbm}}(0,dv)$ of the Sticky Brownian motion (with $\theta= \sqrt{2}\gamma$), namely
\beq
p_t^{\text{sbm}}(0,dv)= \sqrt 2 \gamma e^{2\sqrt 2 \gamma |v|+4\gamma^2t} \erf\(2\gamma \sqrt t +\frac{|v|}{\sqrt{2t}}\)\, dv + \delta_0(dv) e^{4\gamma^2 t}\; \erf\(2\gamma \sqrt t\)
 \eeq
This concludes the proof.
\epr

\br
Using  the expression of  the Laplace transform of  $p_t^{\text{sbm}}(0,dv)$ given in Section 2.4 of \cite{howitt2007stochastic}, it is possible to verify that the Laplace transform of  \eqref{timevarthmRM} (using \eqref{ciao}) coincides with the expression   in Theorem 2.18 of \cite{carinci2017exact}. 
\er

\subsection{Proof of Theorem \ref{MainThm}: Mosco convergence for inclusion dynamics}\label{MCInclusec}

In this section we prove Theorem \ref{MainThm}; the Mosco convergence of the Dirichlet forms associated to the difference process $\{w(t), \: t\ge 0\}$ with infinitesimal generator \eqref{wgen} to the Dirichlet form corresponding to the two-sided sticky Brownian motion. We take the limit in the sticky regime introduced earlier in Section \ref{sect:cond}. In this regime the corresponding scaled difference process is given by:
\be
w_N(t):= \frac 1 N \; w\left(\frac {N^3 \gamma}{\sqrt{2}}\, t\right) \qquad \text{with inclusion-parameter} \quad  k_N:= \frac 1 {\sqrt 2 \gamma  N} \nn
\ee
with infinitesimal generator 
\be\label{wgenN}
    (L_N f) (w) = \frac{N^3 \gamma}{\sqrt{2}} \sum_{r \in A_N } 2 p_N(r) \left(  \frac{1}{ \sqrt{2}N \gamma} + \1_{r=-w} \right) \left[ f(w+r)-f(w)\right]
\ee
 for $w \in \frac{1}{N} \Z$, with
\be\label{AN}
p_N(r):=p(Nr) \qquad \text{and} \qquad
A_N := \frac{1}{N} \{ -R,-R +1, \dots,R-1, R \} \setminus \{0\}
\ee
Notice that by Proposition \ref{reversinurho} the difference processes are reversible with respect to the measures $\nu_{\gamma,N} $ given by
\be\label{nu_distN}
\nu_{\gamma,N} =   \mu_N + \sqrt{2} \gamma \delta_0 
\ee
and by \eqref{DirGenRel} the corresponding sequence of Dirichlet forms is given by
\beq\label{DirwN}
\caE_N (f) =-  \sum_{w \in \Z/N} f(w) \sum_{ r \in A_N} 2 p_N(r) \left( \tfrac{N^2}{2} + \tfrac{N^3 \gamma}{\sqrt{2}} \1_{r=-w}  \right)  (f(w+r)-f(w)) \ \nu_{\gamma,N}  (w)\nn
\eeq

\br\label{RMmeasure}
Notice that  the choice of the reversible measures $\nu_{\gamma,N}$ determines  the sequence  of approximating  Hilbert spaces given by $H_N^{\text{sip}} :=L^2( \frac{1}{N}\Z, \nu_{\gamma,N})$, $N\in \N$. Here for $f, g \in H_N^{\text{sip}}$ their inner product  is given by
\beq\label{L2NSip}
\langle f, g \rangle_{H_N^{\text{sip}}} &=& \sum_{w \in \Z/N} f(w) g(w) \ \nu_{\gamma,N}  (w)
=  \langle f, g \rangle_{H_N^{\text{rw}}} + \sqrt{2} \gamma f(0) g(0)
\eeq
where
\be
\langle f, g \rangle_{H_N^{\text{rw}}} = \frac{1}{N} \sum_{w \in \Z/N} f(w) g(w) \nn
\ee
is the inner product of Section \ref{warmsec}.
\er

\subsubsection{Convergence of Hilbert spaces}

As we already mentioned in Remark \ref{RMmeasure}, by choosing the reversible measures $\nu_{\gamma,N}$ we have determined the convergent sequence of Hilbert spaces and, as a consequence, we have also set the limiting Hilbert space $H^{\text{sbm}} $ to be  $L^2( \R, \bar\nu)$ with $\bar \nu$ as in \eqref{DomSBMchi}.
Notice that from the regularity of this measure, by Theorem 13.21 in \cite{hewitt1975real} and standard arguments we know that the set $C_k^{\infty} (\R)$ of smooth compactly supported test functions is dense in $L^2( \R, \bar\nu)$. Moreover the set
\be
C^0 ( \R \setminus \{ 0\}) := \{ f + \lambda \1_{ \{0 \}} : f \in  C_k^{\infty} (\R) , \lambda \in \R  \}
\ee
denoting  the set of all continuous functions on $\R \setminus \{ 0 \}$ with finite value at $0$, is also dense in $L^2( \R, \bar\nu)$.\\

\noindent
We have to define the right "embedding" operators $\{ \Phi_N \}_{N \geq 1}$, cf. Definition \ref{HilConv}  , to not only guarantee convergence of Hilbert spaces $H_N \to H$ , but Mosco convergence as well. We define these operators as follows:
\be\label{PHINSIP}
\{\Phi_N: C^0 ( \R \setminus \{ 0\}) \to H^{\text{sip}}_N \}_N \qquad \text{defined by}\qquad    \Phi_N f = f \mid_{\Z/N}.
\ee

\bp\label{HilSIPcon}
The sequence of spaces $H_N^{\text{sip}}=L^2( \frac{1}{N}\Z, \nu_{\gamma,N}  )$, $N\in\mathbb N$, converges, in the sense of Definition \ref{HilConv}, to the space $H^{\text{sbm}} =L^2( \R, \bar\nu)$.
\ep

\bpr
The statement follows from the definition of $\{ \Phi_N \}_{N \geq 1}$.
\epr

\subsubsection{Mosco I}\label{MoscoISIPproof}

We will divide our task in two steps. First, we will compare the inclusion Dirichlet form with a random walk Dirichlet form  and show that the first one dominates the second.  We will later use this bound and the fact that the random walk Dirichlet form  satisfies Mosco I, to prove that Mosco I also holds for the case of inclusion particles.\\

\noindent
We consider a random walk on $\Z$ with jump range $A=[-R,R]\cap \Z/\{0\}$. We call again $\{v(t), \;t\ge 0\}$ this process, as in the case of nearest-neighbor RW (that is a special case of this process corresponding to the choice $R=1$). More generally,  in this section we will use the same notation that has been used in Section \ref{warmsec} for the case $R=1$, thus we denote by $L^{rw}$ the infinitesimal  generator:
\be
    (L^{rw}f) (v) = \sum_{r \in A } p(r) \left[ f(v+r)-f(v)\right], \qquad v \in \mathbb Z
\ee
Hence, in the diffusive scaling, the $N$-infinitesimal generator is given by:
\be\label{FRRWDGEN}
 \Delta_N g (v) = N^2 \sum_{r \in A_N^+ }  p_N(r)   \left[g(v+r)-2 g(v) + g(v-r)\right], \qquad v \in \tfrac{\mathbb Z}N
\ee
where $A_N^+ := \{ \rvert r \lvert : r \in A_N   \}$
i.e. the generator of the process $v_N(t):=\tfrac 1 N v(N^2t)$, $t\ge 0$, and  denote by $(\caR_N,D(\caR_N))$ the associated Dirichlet form.\\

\subsubsection*{Comparing RW and SIP Dirichlet forms}
The key idea to prove Mosco I is to transfer the difficulties of the SIP nature to independent random walkers. This is done by means of the following observation:
\bp\label{SIPVSIRW}
For any $f_N \in H_N^{\text{sip}} $ we have
\be
\caE_N (f_N ) \geq  \caR_N (f_N )
\ee
\ep
\bpr
 Rearranging \eqref{DirwN} and using the symmetry of $ p(\cdot)$ allows us to write:
\be
\caE_N (f_N ) -  \caR_N (f_N ) =\frac{N^2}{\sqrt{2}} \,  \gamma \sum_{r \in A_N} 2 p_N(r) (f_N(r)-f_N(0))^2
\ee
and the result  follows from the fact that the RHS of this identity is nonnegative.
\epr

%
%
%

\subsubsection*{Strong and weak convergence in $H_N^{\text{rw}}$ and $H_N^{\text{sip}}$ compared}

\bp\label{strongconvIndzeroRW}
The sequence  $\{ h_N =\1_{\{0\}} \}_{ N \geq 1}$,  with $h_N \in H_N^{\text{rw}}$, converges strongly to $h=0 \in H^{\text{bm}}$ with respect to $H_N^{\text{rw}}$-Hilbert convergence.
\ep

\bpr
In the language of Definition \ref{strongcon} we set $\tilde{h}_M\equiv 0 $. With this choice we immediately have
\be
\| \hat{h}_M -h \|_{H^{\text{bm}}} = 0
\qquad \text{and} \qquad
\| \Phi_N \hat{h}_M -h_N \|_{H_N^{\text{rw}}}^2 = \tfrac{1}{N}
\ee
which concludes the proof.
\epr

\bp\label{strongconvIndzero}
The sequence  $\{ h_N =\1_{\{0\}} \}_{ N \geq 1}$,  with $h_N \in H_N^{\text{sip}}$, converges strongly to $h=\1_{\{0\}} \in H^{\text{sbm}}$ with respect to $H_N^{\text{sip}}$-Hilbert convergence.
\ep

\bpr
In the language of Definition \ref{strongcon} we set $\tilde{h}_M\equiv \1_{\{0\}} $. With this choice we immediately have
\be
\| \hat{h}_M -h \|_{H^{\text{sbm}}} = 0
\qquad \text{and} \qquad
\| \Phi_N \hat{h}_M -h_N \|_{H_N^{\text{sip}}} = 0
\ee
which concludes the proof.
\epr

\noindent
A consequence of Proposition \ref{strongconvIndzero} is that any sequence weakly convergent, with respect to $H_N^{\text{sip}}$-Hilbert convergence, converges also at zero.

\bp\label{WeakconatzeroSIP}
Let $\{ f_N \}_{N \geq 1}$ in $\{ H_N^{\text{sip}} \}_{N \geq 1}$ be a sequence converging weakly  to $f \in H^{\text{sbm}}$ with respect to $H_N^{\text{sip}}$-Hilbert convergence,  then $\lim_{N \to \infty} f_N(0) = f(0)$.
\ep

\bpr
By Proposition \ref{strongconvIndzero} we know that $\{ h_N =\1_{\{0\}} \}_{ N \geq 1}$ converges strongly to $h=\1_{\{0\}}$ with respect to $H_N^{\text{sip}}$-Hilbert convergence. This, together with the fact that $\{ f_N \}_{N \geq 1}$ converges weakly, implies:
\be\label{weakfNHN}
\lim_{N\to \infty} \langle f_N, h_N  \rangle_{H_N^{\text{sip}}} = \  \langle f, h  \rangle_{H^{\text{sbm}}}= \sqrt{2} \gamma f(0)
\ee
but by \eqref{L2NSip}
\be
\langle f_N, h_N \rangle_{H_N^{\text{sip}}} =  ( \tfrac{1}{N}  + \sqrt{2} \gamma ) f_N(0)
\ee
which, together with \eqref{weakfNHN}, implies the statement.
\epr

\noindent
To further contrast the two notions of convergence, Proposition \ref{strongconvIndzeroRW} has a weaker implication
\bp\label{WeakconatzeroRW}
Let $\{ g_N \}_{N \geq 1}$ in $\{ H_N^{\text{rw}} \}_{N \geq 1}$ be a sequence converging weakly  to $g \in H^{\text{bm}}$ with respect to $H_N^{\text{rw}}$-Hilbert convergence, then $\lim_{N \to \infty} \tfrac{1}{N} g_N(0) = 0$.
\ep

\bpr
By Proposition \ref{strongconvIndzeroRW} we know that $\{ h_N =\1_{\{0\}} \}_{ N \geq 1}$ converges strongly to $h=0$ with respect to $H_N^{\text{rw}}$-Hilbert convergence. This, together with the fact that $\{ g_N \}_{N \geq 1}$ converges weakly, implies:
\be\label{weakgNHN}
\lim_{N\to \infty}  \langle g_N, h_N  \rangle_{H_N^{\text{rw}}} =  0
\ee
but we know
\be
\left \langle g_N, h_N \right \rangle_{H_N^{\text{rw}}} =  \frac{1}{N}   g_N(0)
\ee
which together with \eqref{weakgNHN} concludes the proof.
\epr

\subsubsection*{From $H_N^{\text{rw}}$ strong convergence to $H_N^{\text{sip}}$ strong convergence}

\bp\label{strongconequi}
Let $\{ g_N \}_{N \geq 1}$ in $\{ H_N^{\text{rw}}  \}_{N \geq 1}$ be a sequence converging strongly  to $g \in H^{\text{bm}}$ with respect to $H_N^{\text{rw}}$-Hilbert convergence. For all $N\geq 1$ define the sequence
\be\label{gNhat}
\hat{g}_N = g_N  - g_N(0) \1_{\{0\}}
\ee
Then $\{ \hat{g}_N \}_{N \geq 0}$ also converges  strongly  with respect to $H_N^{\text{sip}}$-Hilbert convergence to $\hat{g}$ given by:
\be\label{ghat}
\hat{g} = g  - g(0) \1_{\{0\}}
\ee
\ep

\bpr
From the strong convergence in the $H_N^{\text{rw}} $-Hilbert convergence sense, we know that there exists a sequence $\tilde{g}_M \in C_k^{\infty} (\R)$ such that
\be
\lim_{M\to \infty} \| \tilde{g}_M -g \|_{H^{\text{bm}}} = 0
\ee
and
\be\label{suplimgm}
\lim_{M \to \infty} \limsup_{N \to \infty} \| \Phi_N \tilde{g}_M -g_N \|_{H_N^{\text{rw}}} = 0
\ee
\noindent
for each $M$ we define the function $\hat{g}_M $ given by
\be\label{gmhat}
\hat{g}_M = \tilde{g}_M - \tilde{g}_M(0) \1_{\{0\} }\nn
\ee
Notice that:
\be
\| \hat{g}_M \|_{H^{\text{sbm}}}^2 = \| \tilde{g}_M \|_{H^{\text{bm}}}^2   < \infty
\ee
and hence we have $\hat{g}_M$ belongs to both $C^0 ( \R \setminus \{ 0\})$ and $H^{\text{sbm}}$. \\

\noindent
As before, we have the relation:
\beq
\| \hat{g}_M - \hat{g} \|_{H^{\text{sbm}}}^2 = \| \hat{g}_M - \hat{g} \|_{H^{\text{bm}}}^2 + \sqrt{2}\gamma( \hat{g}_M(0) - \hat{g}(0)  )^2 
= \| \tilde{g}_M - g \|_{H^{\text{bm}}}^2
\eeq
which shows that indeed we have
\be\label{equivstronprop1}
\lim_{M\to \infty} \| \hat{g}_M -\hat{g} \|_{H^{\text{sbm}}}^2 = 0
\ee
For the second requirement of strong convergence we can estimate as follows
\beq
\| \Phi_N \hat{g}_M -\hat{g}_N \|_{H_N^{\text{sip}}}^2 &=& \frac{1}{N} \sum_{\substack{ x \in \tfrac{1}{N} \Z \\ x \neq 0} } ( \Phi_N \tilde{g}_M(x) -g_N(x))^2 \nn
\eeq
relation \eqref{suplimgm} allows us to see that the RHS of the equality above vanish. This, together with \eqref{equivstronprop1} concludes the proof of the Proposition.
\epr

\subsubsection*{From $H_N^{\text{sip}}$  weak convergence to $H_N^{\text{rw}}$ weak convergence}

The following proposition says that with respect to weak convergence the implication comes in the opposite direction

\bp\label{Weakconequi}
Let $\{ f_N \}_{N \geq 1}$ in $\{ H_N^{\text{sip}} \}_{N \geq 1}$ be a sequence converging weakly  to $f \in H^{\text{sbm}}$ with respect to $H_N^{\text{sip}}$-Hilbert convergence. Then it also converges  weakly  with respect to $H_N^{\text{rw}}$-Hilbert convergence.
\ep

\bpr
Let $\{ f_N \}_{N \geq 0}$ in $\{ H_N^{\text{sip}} \}_{N \geq 0}$ be as in the Proposition. In order to show that it also converges weakly with respect to $H_N^{\text{rw}}$-Hilbert convergence we need to show that for any sequence $\{ g_N \}_{N \geq 0}$ in $\{ H_N^{\text{rw}}  \}_{N \geq 0}$   converging strongly to some $g \in H^{\text{bm}}$ we have
\be
\lim_{N\to \infty}  \langle f_N, g_N \rangle_{H_N^{\text{rw}}} =  \langle f, g \rangle_{H^{\text{bm}}}
\ee

\noindent
Consider such a sequence $\{ g_N \}_{N \geq 0}$, by Proposition \ref{strongconequi} we know that the sequence $\{ \hat{g}_N \}_{N \geq 1}$ also converges stronlgy with respect to $H_N^{\text{sip}}$-Hilbert convergence to $\hat{g}$ defined as in \eqref{ghat}. Then we have:
\beq
\lim_{N\to \infty}  \langle f_N, \hat{g}_N  \rangle_{H_N^{\text{sip}}} &=&  \langle f, \hat{g}  \rangle_{H^{\text{sbm}}} =  \langle f, g  \rangle_{H^{\text{bm}}}
\eeq
which can be re-written as:
\beq
\lim_{N\to \infty} \langle f_N, g_N  \rangle_{H_N^{\text{rw}}} - \frac{1}{N} f_N(0) g_N(0) &=&   \langle f, g \rangle_{H^{\text{bm}}}
\eeq
and together with Propositions \ref{WeakconatzeroSIP} and \ref{WeakconatzeroRW}  implies that:
\beq
\lim_{N\to \infty}  \langle f_N, g_N  \rangle_{H_N^{\text{rw}}}  &=&   \langle f, g  \rangle_{H^{\text{bm}}}
\eeq
and the proof is done.
\epr

\subsubsection*{Conclusion of proof of Mosco I}

In order to see that condition Mosco I is satisfied, we combine Proposition \ref{SIPVSIRW}, Proposition \ref{Weakconequi} and the Mosco convergence of Random Walkers to Brownian motion to obtain that for all $f \in H^{\text{sbm}}$, and all $f_N \in H_N^{\text{sip}}$ converging weakly to $f$ we have
\be
\liminf_{N \to \infty} \caE_N (f_N ) \geq \liminf_{N \to \infty} \caR_N (f_N ) \geq \caE_{\text bm}(f) =\caE_{\text sbm}(f) \nn
\ee
where the last equality comes from Remark \ref{remarkEqualDirs}.

\subsubsection{Mosco II}
We are going to prove that Assumption 2 is satisfied. We use the set of compactly supported smooth functions $C_k^\infty (\R)$, which by the regularity of the measure $dx+ \delta_0$ is dense in $H=L^2(dx+ \delta_0)$.

\subsubsection*{The recovering sequence} 

For every $f \in C_k^\infty (\R)$, we need to find a sequence $f_N$ strongly-converging  to $f$ and such that
\be\label{aimMII}
\lim_{ N \to \infty} \caE_N ( f_N ) = \caE (f)
\ee
The obvious choice $f_N = \Phi_N f$ does not work in this case, the reason of this is the emergence in the limit of a non-vanishing term containing $f\myprime(0)$. Nevertheless our candidate is the sequence $\{ \Psi_N f \}_{N \geq 1}$ given by
\begin{equation}\label{PSIN}
    (\Psi_N f ) (i) =
\left\{
    \begin{array}{ll}
        f(i) & i \in \tfrac{1}{N}\Z\setminus A_N \\
        f(0) & \text{otherwise}
    \end{array}
\right.
\qquad \quad \text{for any} \: f \in C_k^\infty (\R)
\end{equation}
where $A_n$ is as in \eqref{AN}. 

\br
The sequence $\{ \Psi_N f \}_{N \geq 1}$ is chosen in such a way that the SIP part of the Dirichlet form, i.e. $\caE^{\text{sip}} - \caR$, vanishes at $\Psi_N f$ for all $N$. See below for the details.
\er

\noindent
Our goal is to show that the sequence $\{ \Psi_N f \}_{N \geq 1}$ indeed satisfies \eqref{aimMII}. First of all we need to show that $\Psi_N f \to f$ strongly.

\bp\label{RecoverySeq}
For all $f \in C_k^\infty (\R) \subset L^2(dx+ \delta_0)$, the sequence $\{ \Psi_N f \}_{N \geq 1}$ in $H_N^{\text{sip}}$ strongly-converges to $f$ w.r.t. the $H_N^{\text{sip}}$-Hilbert space convergence given.
\ep

\bpr
In the language of Definition \ref{strongcon} we set $\tilde{f}_M\equiv f $. Hence the first condition is trivially satisfied:
\be
\lim_{M \to \infty} \| \tilde{f}_M -f \|_{H^{\text{sbm}}}  = 0
\ee
Moreover
\beq
&&\lim_{M \to \infty} \limsup_{N \to \infty} \| \Phi_N \tilde{f}_M - \Psi_N f \|_{H_N^{\text{sip}}}^2 = \limsup_{N \to \infty} \| \Phi_N f - \Psi_N f \|_{H_N^{\text{sip}}}^2 \nn \\
&&= \limsup_{N \to \infty} \sum_{i \in \frac{1}{N}\Z} ( \Phi_N f(i) -\Psi_N f (i) )^2 \nu_{\gamma,N}  (i) = \limsup_{N \to \infty}  \frac{1}{N} \sum_{ i \in A_N}( f(i)-f(0))^2  = 0 \nn
\eeq
where we used the boundedness of $f$ and the fact that the cardinality of the set $A_N$ is finite and does not depend on $N$.
\epr

\subsubsection*{Preliminary simplifications}
To continue the proof of \eqref{aimMII}, the first thing to notice is that the Dirichlet form $\caE_N$ evaluated in $\Psi_N f$ can be substantially simplified:
\beq\label{FRMIIstep1}
\caE_N (\Psi_N f) &=&  -\sum_{i \in \frac{1}{N}\Z} \Psi_N f (i) \sum_{ r \in A_N} 2 p_N(r) \left( \tfrac{N^2}{2} + \tfrac{N^3 \gamma}{\sqrt{2}} \1_{r=-i}  \right)  (\Psi_N f(i+r)- \Psi_N f(i)) \nu_{\gamma,N}  (i) \nn \\
&=& - \sum_{i \in \frac{1}{N}\Z}\Psi_N f (i)  \sum_{ r \in A_N}  p_N(r)  N^2   (\Psi_N f(i+r)- \Psi_N f(i)) \nu_{\gamma,N}(i)  \\
&&-  \sum_{i\in \frac{1}{N}\Z} \Psi_N f (i)  \sum_{ r \in A_N} 2 p_N(r) \left( \tfrac{N^3 \gamma}{\sqrt{2}} \1_{r=-i}  \right)  (\Psi_N f(i+r)- \Psi_N f(i)) \nu_{\gamma,N}  (i)\nn
\eeq
where, from the observation that for $i = -r$ and $r \in A_N$, via \eqref{PSIN} we get
\be
(\Psi_N f(i+r)- \Psi_N f(i)) = 0
\ee
the whole second sum in \eqref{FRMIIstep1} vanishes. Then by \eqref{nu_distN} we are left with
\beq\label{FRMIIstep2}
\caE_N (\Psi_N f)
= &&-N\sum_{ r \in A_N}  p_N(r)  \sum_{i \in \frac{1}{N}\Z}\Psi_N f (i) (\Psi_N f(i+r)- \Psi_N f(i)) \nn \\
&&- \sqrt{2} \gamma  N^2 \sum_{ r \in A_N}  p_N(r) \Psi_N f (0) (\Psi_N f(r)- \Psi_N f(0))
\eeq
 we have again that
$
(\Psi_N f(r)- \Psi_N f(0)) = 0
$
for $r\in A_N$, then
 our Dirichlet form becomes
\beq\label{FRMIIstep3}
&&\caE_N (\Psi_N f) = -N\sum_{ r \in A_N}  p_N(r)  \sum_{i \in \frac{1}{N}\Z}   \Psi_N f (i) (\Psi_N f(i+r)- \Psi_N f(i)) \nn
\eeq
that we split again as follows
\beq\label{FRMIIstep4}
&&\caE_N (\Psi_N f) = -N\sum_{ r \in A_N}  p_N(r) \sum_{\substack{i \in \frac{1}{N}\Z \setminus A_N }}   \Psi_N f (i) (\Psi_N f(i+r)- \Psi_N f(i))-S_N \nn \\
&&\text{with}\qquad S_N=N \sum_{ r \in A_N}  p_N(r)  \sum_{ i \in A_N }  \Psi_N f (i) (\Psi_N f(i+r)- \Psi_N f(i))
\eeq
\subsubsection*{The correct limit}
First we show that $S_N$ vanishes as $N \to \infty$. For $i \in A_N$, we define the sets
\beq
A_N^i:= A_N -i \qquad
\text{and} \qquad
A_N^+ = \{ \rvert r \lvert : r \in A_N   \}
\eeq
notice that for  $ r \in A_N^i$ we have
$(\Psi_N f(i+r)- \Psi_N f(0)) = 0$
and hence
\beq\label{secondterm1}
S_N &&=  N\sum_{ i \in A_N } \sum_{ r \in A_N \setminus A_N^i} p_N(r)  f (0)     ( f(i+r)-  f (0)) \nn \\
&&=  N\sum_{ i \in A_N^+ } \sum_{ r \in A_N \setminus A_N^i} p_N(r)  f (0)    ( f(i+r)- 2 f(0) + f(-i-r) ) \nn
\eeq
where we used the symmetry of $p(\cdot)$ and the fact that  $r \in A_N \setminus A_N^i$ if and only if $-r \in A_N \setminus A_N^{-i}$. We conclude that $S_N$ vanishes by recalling that by a Taylor expansion the factor $( f(i+r)- 2 f(0) + f(-i-r) )$ is of order $N^{-2}$.\\

\noindent
For what concerns the remaining term in \eqref{FRMIIstep4}, we notice that, exploiting  the   symmetry of the transition function $p(\cdot)$, we can re-arrange it into
\beq
\caE_N (\Psi_N f) {+S_N}
= - N\sum_{ r \in A_N^+}  p_N(r) \sum_{\substack{i \in \frac{1}{N}\Z \setminus A_N  }}  \Psi_N f (i) (\Psi_N f(i+r)-2 \Psi_N f(i)+\Psi_N f(i-r))\nn
\eeq
Let us define the following set $B_N =\frac{1}{N} \{ -2R,-2R +1, \dots,2R-1, 2R \}$
and split the sum above as follows
\beq\label{FRMIIstep6}
&&\caE_N (\Psi_N f) {+S_N} = -N\sum_{ r \in A_N^+}  p_N(r) \sum_{\substack{i \in \frac{1}{N}\Z \setminus B_N }}  \Psi_N f (i) (\Psi_N f(i+r)-2 \Psi_N f(i)+\Psi_N f(i-r)) \nn \\
&&\hskip2cm - N\sum_{ r \in A_N^+}  p_N(r)  \sum_{\substack{i \in B_N \setminus A_N  }}  \Psi_N f (i) (\Psi_N f(i+r)-2 \Psi_N f(i)+\Psi_N f(i-r)) 
\eeq
The above splitting allows to isolate  the first term for which we have no issues of the kind $\Psi_N f( i +r ) = f(0)$ and hence no complications when taylor expanding around the points $i \in \frac{1}{N}\Z$. \\

\noindent
We now show  that the second term in the RHS  of \eqref{FRMIIstep6} vanishes as $N$ goes to infinity: \\

\noindent
Take a positive $ i \in B_N \setminus A_N$, then for $ r \in A_N^i$,
$\Psi_N f (i+r) = f(0)$.

\br
Notice that, for $-i  \in B_N \setminus A_N$, the set $A_N^{-i}= A_N^{i}$  is such that
\be
\Psi_N f (-i + r) = f(0)\qquad \text{for all} \quad  r \in A_N^{i}.
\ee

\er
\br
We will omit the analysis for $ r\notin A_N^i$ because for those terms we can Taylor expand $f$ around the point $i$ and show that the factors containing the discrete Laplacian are of order $N^{-2}$.
\er
\noindent
We now consider the contribution that each pair  $(i,-i)$ gives to the second sum in the RHS of \eqref{FRMIIstep6}. Let $ i \in (B_N \setminus A_N )^+$, then
\beq\label{contriplus}
C_N(i) && :=N \sum_{ r \in A_N^i}  p_N(r)   \Psi_N f (i) \left[ \Psi_N f(i+r)-2 \Psi_N f(i)+\Psi_N f(i-r) \right] \nn \\
&&=N \sum_{ r \in A_N^i}  p_N(r)    f (i) \left[  f(i+r)-2 f(i)+ f(0) \right]
\eeq

\noindent
Taylor expanding around zero the terms inside the square brackets in the RHS of \eqref{contriplus} gives
\beq\label{contriplus2}
C_N(i) &&= \sum_{ r \in A_N^i}  p_N(r)   f (i) f\myprime(0) \left[ r -i  \right] + O(1/N) \nn
\eeq

\noindent
Analogously, for the contribution $C_N(-i)$ we obtain
\beq\label{contriminus2}
C_N(-i) &&= \sum_{ r \in A_N^i}  p_N(r)   f (-i) f\myprime(0) \left[ i- r \right] + O(1/N) \nn
\eeq
summing both contributions over all $i >0$ we obtain
\beq\label{contrissum}
&&\sum_{\substack{i \in (B_N \setminus A_N )^+ }} C_N(i) + C_N(-i)\\
&&= \sum_{\substack{i \in (B_N \setminus A_N)^+  }} \sum_{ r \in A_N^i}  p_N(r)  f\myprime(0)  \left( r -i  \right) \left[ f (i) - f (-i)  \right] + O(1/N) = O(1/N)\nn
\eeq
where we used that the cardinality of the sets $A_N^i$ and $(B_N \setminus A_N )^+$ does not depend on $N$. Then we can write
\beq\label{FRMIIstep7}
&&\caE_N (\Psi_N f)= - \frac{1}{N} \sum_{ r \in A_N^+}  p_N(r)\sum_{\substack{i \in \frac{1}{N}\Z \setminus B_N }}  N^2  f (i) ( f(i+r)-2  f(i)+ f(i-r)) + O(1/N) \nn
\eeq
which indeed by a Taylor expansion gives the limit
\beq
\lim_{ N\to \infty } \caE_N (\Psi_N f) = -\frac{\chi}{2}  \int_{\R} f(x) f\mydprime(x) \ dx = \frac{\chi}{2}  \int_{\R} f\myprime(x)^2 \ dx
\eeq
with $\chi = \sum_{r=1}^R  p(r) r^2 $. This concludes the proof of Mosco II. \qed

\newpage

\section{Appendix}

%

\subsection{Sticky Brownian Motion and its Dirichlet form}

In this Appendix we provide some background material on the two sided sticky Brownian motion in the context of Dirichlet forms . Namely, by means of an example we apply the machinery of Dirichlet forms to the theory of stochastic time changes for Markov processes. The example that we will build at the end of this section plays the role of the limiting process for the difference process. In this appendix will mostly follow the approach presented in Chapter 5 of \cite{chen2012symmetric}.

\subsubsection{Two sided sticky Brownian motion}

The traditional approach to construct sticky Brownian motion (SBM) on the real line is by means of local times and time changes related to them. Let us say that we are in the one dimensional case and we want to build  Brownian motion sticky at zero. We consider then standard Brownian motion $\{ B_t \}_{ t \geq 0}$ taking values on $\R$ and define its local time at zero by
\be
L_t^0 = \lim_{ \epsilon \to 0} \frac{1}{2\epsilon} \int_0^t \1_{[-\epsilon,\epsilon]} (B_s) ds
\ee
Given this local time and for $\gamma > 0$ we consider the functional
\be\label{acctime}
T_t = t + \gamma L_t^0
\ee
and denote by $\tau$ its generalized inverse, i.e.,
\be\label{tauinv}
\tau(t)= \inf \{ s > 0 :  T_s > t  \}
\ee
then the process given by the time change
\be
B_t^{\text{sbm}} = B_{\tau(t)}
\ee
is what is known in the literature by two sided sticky Brownian motion.
\br
The idea in defining \eqref{acctime} is that we add some \quotes{extra time} at zero and by taking the inverse \eqref{tauinv} via the time change we slow down the new process whenever it is at 0. Notice that the parameter $\gamma$ controls the factor by which we slow down time.
\er
\noindent
As expected, in the context of Dirichlet forms, we can also perfom this kind of stochastic time changes. Our goal for this section is to describe the Dirichlet forms approach to perfom the kind of time changes we are interested in. There are basically two ingredients that we need:
\ben
\item A symmetric Markov process $M_t$ with reversible measure $\mu$ with  support in the state space $E$.
\item A Positive Continous Additive Functional (PCAF) that, in a sense to be seen later, plays the role of the local time.
\een

\br
In the same way that the local time $L_t^0$ implicitely defined the point $\{0\}$ as the \quotes{sticky region}, the PCAF of the second ingredient above will determine a \quotes{sticky region} for our new process.
\er

\noindent
For the sake of completeness let us introduce the precise definition of PCAF's

\bd[PCAF]
A function $A_t (\omega)$ of two variables $t\geq 0$ and $\omega \in \Omega$ is called an additive functional of $M_t$ if there exists $\Lambda \in \caF_\infty$ and a $\mu$-inessential set $N \subset E$ with
\be
P_x (\Lambda) = 1 \quad \text{for } x \in E\setminus N \quad \text{and} \quad \theta_t \Lambda \subset \Lambda \quad \text{ for } t \geq 0
\ee
and the following conditions are satisfied:
\begin{description}
\item[(i)] For each $ t \geq 0$, $A_t \mid_{\Lambda}$ is $\caF_t$-measurable.
\item[(ii)] For any $\omega \in \Lambda$, $A_{\cdot} (\omega)$ is right continuous on $[0,\infty)$ has left limits on $(0,\zeta(\omega))$, $A_0 (\omega)=0$, $\lvert A_t (\omega) \rvert < \infty$ for $t < \zeta(\omega)$, and $A_t (\omega)= A_{\zeta(\omega)} (\omega)$ for all $ t \geq \zeta(\omega)$.
\item[(iii] The additivity property is satisfied, i.e.,
\be
A_{t+s} (\omega) = A_{t} (\omega) + A_{s} (\omega) \text{ for all } t,s \geq 0
\ee
\end{description}
\ed

\noindent
If we denote by $\caA_c^{+}$ the set of all PCAF, it turns out that there exists a one to one correspondence between the set $\caA_c^{+}$ and a special subset of the set of the Borel measures on $E$. Which we now introduce:

\bd[Smooth measures]
Let $\nu$ be a positive measure on $(E, \caB(E))$, $\nu$ is said to be smooth if
\ben
\item It does not charge any $\caE_M$-polar set.
\item There exists a nest $\{ F_k \}_{k \geq 1}$ such that $\nu (F_k) < \infty$ for all $k \geq 1$.
\een
\ed

\br
Notice that all the Dirichlet forms related concepts ( $\caE_M$-capacity  for example ) are in terms of the Dirichlet space $(\caE_M,D(\caE_M))$, which corresponds to the symmetric Markov process $M_t$.
\er
\noindent
We denote by $S(E)$ the set of all smooth measures on $E$. The correspondence we mentioned above is between $\caA_c^{+}$ and $S(E)$. Formally, this correspondence is given by the following result:

\bt[PCAF and Smooth measures]
For $A \in \caA_c^+$ we denote by $\nu_A$ the measure that is in Revuz correspondence with $A$, i.e. the measure that for any $ f \in \caB_{+} (E)$ satisfies:
\be\label{Revuzrel}
\int_E f (x) \nu_A (dx) = \lim_{t \downarrow 0} \frac{1}{t} E_{\mu} [\int_0^t f (M_s) d A_s ]
\ee
then we have the following:
\begin{description}
\item[(i)] For any $A \in \caA_c^+$, $\nu_A \in S(E)$.
\item[(ii)] For any $\nu \in S(E)$, there exists $A \in \caA_c^+$ satisfying $\nu_A = \nu$ uniquely up to $\mu$-equivalence.
\end{description}
\et

\bpr
This is part of Theorem 4.1.1 in \cite{chen2012symmetric} where the proof is included.
\epr

\noindent
It is known that there exists a one to one correspondence between Markov process and Dirichlet forms \cite{fukushima2010dirichlet}. The idea is that given a PCAF $A_t$ we can define a stochastic time changed process given by the generalized inverse of $A_t$ in terms of its corresponding Dirichlet form. More precisely:

\bt\label{DirTchange}
Let $M_t$ be a symmetric Markov process with corresponding Dirichlet space given by $(\caE_M,D(\caE_M))$. Let also $A_t$ be a PCAF whose Revuz measure $\nu_A$ has full quasi support. Denote by $\tilde{M}_t$ the time changed process given by the generalized inverse of $A_t$. Then we have that its corresponding  Dirichlet space $(\caE_{\tilde{M}},D(\caE_{\tilde{M}}))$ is given by
\be
\caE_{\tilde{M}} (f,g) = \caE_{M} (f,g)  \quad  \text{and} \quad D(\caE_{\tilde{M}}) = D(\caE_M) \cap L^2 (E,\nu_A)
\ee
\et

\bpr
This theorem is just a specialization of Theorem 5.2.2 in \cite{chen2012symmetric}. Where the time changed form is given by
\be
\caE_{\tilde{M}} (f,g) = \caE_{M} (H_F f , H_F g)
\ee
The specialization consists in the fact that the Revuz measure $\nu_A$ has full quasi support, i.e.,
\be
H_F h( x) = \E_x [ h(M_{\sigma_F} ) ; \sigma_f < \infty ] = h(x)
\ee
where $F$ is the support of $\nu_A$ and $\sigma_F$ is its hitting time. We refer the reader to page 176 of the same reference if more details are needed.
\epr

\noindent
Under this setting, it becomes then easier to characterize the time changed of Brownian motion given by the inverse of the functional $T_t$ defined in \eqref{acctime}. The idea is that under the setting given by one dimensional Brownian motion on the reals. We know that the process $\{ B_t \}_{t \geq 0}$ is reversible with respect to the Lebesgue measure $dx$. On the first hand, the Lebesgue measure $dx$ is in Revuz correspondence with the trivial PCAF $A_t^1=t$. On the other hand the following computation shows the Revuz correspondence between the PCAF $L_t^0$ and the Dirac measure at zero $\delta_0$:
\small
\begin{align}\label{dirac0proof}
\lim_{t \downarrow 0} \frac{1}{t} \E_{dx} [\int_0^t f (B_s) d L_s^0 ]
&= \lim_{t \downarrow 0} \frac{1}{t} \E_{dx} [\int_0^t f (B_s) \lim_{ \epsilon \downarrow 0} \frac{1}{2\epsilon} \1_{[-\epsilon,\epsilon]} (B_s) \ ds] \nn \\
&= \lim_{t \downarrow 0} \lim_{ \epsilon \downarrow 0} \frac{1}{t} \frac{1}{2\epsilon}  \int_0^t \int_{\R} \E_{B_0} [ f (B_s+x)  \1_{[-\epsilon,\epsilon]} (B_s+x) ] \ dx ds \nn \\ \nn \\
&= \lim_{t \downarrow 0} \lim_{ \epsilon \downarrow 0} \frac{1}{t} \frac{1}{2\epsilon}  \int_0^t \int_{\R^2}  f (y+x)  \1_{ [-\epsilon,\epsilon]}(y+x) \frac{e^{\tfrac{-y^2}{2s}}}{\sqrt{2 \pi s}} \  dy dx ds \nn \\ \nn \\
&= \lim_{t \downarrow 0} \lim_{ \epsilon \downarrow 0} \frac{1}{t} \frac{1}{2\epsilon}  \int_0^t \int_{\R} \int_{-\epsilon}^{\epsilon}  f (z) \frac{e^{\tfrac{-(z-x)^2}{2s}}}{\sqrt{2 \pi s}} \  dz dx ds \nn \\ \nn \\
&=  \lim_{ \epsilon \downarrow 0}  \frac{1}{2\epsilon}    \int_{-\epsilon}^{\epsilon}  f (z) \  dz  \nn \\ \nn \\
&= f(0) = \int f(x) \delta_0 (dx)
\end{align}
\normalsize
\noindent
Then the measure $\nu = dx + \gamma \delta_0$ is in Revuz correspondance with the PCAF $T_t$ and hence by Theorem \ref{DirTchange} the Dirichlet form for one dimensional Sticky Brownian motion $\{ B_t^{\text{sbm}} \}_{ t \geq 0}$ is given by:
\be\label{timchanSBM}
\caE_{B^{\text{sbm}}} (f,g) = \caE_{B} (f,g)  \quad  \text{and} \quad  D(\caE_{B^{\text{sbm}}}) = D(\caE_{B}) \cap L^2 (\R ,dx + \gamma \delta_0)
\ee

\subsubsection{Domain of the infinitesimal generator}

With the objective to obtain a description of the generator of the time changed process that we have just built.  In this section we will make use of the correspondence between Dirichlet forms and Markov generators. Let us then expand a bit on what we mentioned before equation \eqref{DirGenRel}; this is how the two directions of the correspondence are actually given:
\begin{description}
\item[(a) From forms $\caE$ to generators $L$:] The correspondence is defined by
\beq\label{correspa}
D(L) \subset D(\caE), \quad \caE(f,g)= -<Lf, g> \quad \forall f \in D(L), \, g \in D(\caE)
\eeq
\item[(b) From generators $L$ to forms $\caE$:] In this case the correspondence is given by
\beq\label{correspb}
D(\caE)=D(\sqrt{-L}), \quad \caE(f,g)= <\sqrt{-L}f, \sqrt{-L}g> \quad \forall f, g \in D(\caE)
\eeq
\end{description}
\noindent
We can think of these relations as the first and second representation theorems for Dirichlet forms in the spirit of Kato \cite{kato2013perturbation} for sesquilinear forms. For the particular case of Dirichlet forms, more details and the connection to semigroups and resolvents, can be found on the Appendix of \cite{chen2012symmetric}. \\

\br
Please notice that the time changed process behaves like Brownian motion on the set $\R \setminus\{ 0\}$ and differently (sticky behavior) when it visits $0$. Therefore we expect the new generator $L_{B^{\text{smb}}}$ to be the same Laplace operator in the region$\R \setminus\{ 0\}$ i.e.
\beq\label{comparingLDelta}
L_{B^{\text{sbm}}} f(x) = \Delta f (x) \quad \forall x \in \R^2
\eeq
and some additional restrictions at the point zero.
\er

\noindent
The idea is to assume that the generaor $L_{B^{\text{sbm}}}$ is just the Laplacian at all points, and by using the properties of the time changed process determine additional constrains at zero.\\

\noindent
For  $f \in D(\caE_{B^{\text{sbm}}})$. Thanks to \eqref{correspb} we can re-write \eqref{timchanSBM} in  terms of $L_{B^{\text{sbm}}}$ in the following way:
\be\label{findgen}
\caE_{B^s} (f,g) =  \int_{\R \setminus\{ 0\}}  g\myprime(x)  f\myprime(x) dx
\ee
for all $g \in D(\caE_{B^{\text{sbm}}})$. \\

\noindent
On the other hand
\be\label{findgen1}
\caE_{B^s} (f,g) = -\int_{\R \setminus\{ 0\}} g(x) f\mydprime(x) dx   - \gamma  g(0) f\mydprime(0)
\ee
where we took $f$ as a member of $D(L_{B^{\text{sbm}}})$ and used \eqref{correspa}.\\

\noindent
Let us split the first therm on the r.h.s. of \eqref{findgen1} in two regions:
\be\label{findgensplit}
\int_{\R \setminus\{0\}} g(x)  f\mydprime(x) dx  = \int_{x > 0} g(x)  f\mydprime(x) dx  + \int_{x < 0} g(x)  f\mydprime(x) dx
\ee

\noindent
Integrating by parts in the first integral of the r.h.s. of \eqref{findgensplit} we obtain:
\be\label{findgen2}
\int_{x > 0} g(x)  f\mydprime(x) dx = -g(0) f\myprime(0+) -\int_{x > 0}  g\myprime(x) \, f\myprime(x) dx
\ee
where
\be
f^{\prime}(0+) = \lim_{ h \downarrow  0} \frac{ f(h)-f(0)}{h}
\ee

\noindent
Similarly we obtain:
\be\label{findgen3}
\int_{x < 0} g(x) f\mydprime(x) dx = g(0) f\myprime(0-) -\int_{x < 0} g\myprime(x) \, f\myprime(x) dx
\ee
therefore, for every $g \in D(\caE_{B^s})$  we obtain:
\be
g(0) \( \gamma \Delta f(0) - f^{\prime}(0+) + f^{\prime}(0-) \) = 0
\ee
which gives
\be\label{BCSBM2D}
 \gamma f^{\prime\prime}(0) =f^{\prime}(0+) - f^{\prime}(0-)
\ee
for every  $f \in D(L_{B^{\text{sbm}}})$
\br
Notice that condition \eqref{BCSBM2D} coindices with what we would expect from the conditions given for two sided sticky Brownian motion. See for instance Appendix 1 in \cite{borodin2012handbook}.
\er

\subsection{Mosco convergence for the Random Walk}\label{warmsec}

In this section, we consider  the difference process for the position-coordinates of two particles performing nearest-neighbor symmetric independent random walks. This process, that we denote by $\{v(t), t\ge 0\}$, is itself a random walk in $\Z$ for which  convergence to the  standard Brownian motion in the diffusive time-scales is well-known. By convergence we mean convergence of generators. In this section we will prove  Mosco convergence of  Dirichlet forms of $v(t)$.

As we can see in Section \ref{MCInclusec}, the proof of Mosco-convergence for inclusion walkers  strongly relies on the result for independent walkers (in particular for the proof of Mosco I). The choice of considering the independent dynamics case has the purpose to
exemplifying the use of the Dirichlet approach in a setting simpler than the one of inclusion dynamics.\\

\noindent
The  generator of $\{v(t), \; t\ge 0\}$ is given by the discrete Laplacian $\Delta_1$:
\beq\label{LdifftIRW}
L^{\text{rw}}f(v) = \Delta_1 f(v)= f(v+1) - 2f(v)+f(v-1), \qquad v\in\mathbb Z.
\eeq
that is simply  the generator of a random walk in $\Z$. Speeding up time by a factor $ N^2 $ and scaling the mesh between the lattice sites by a factor $\tfrac 1 N$ we obtain that the generator of this scaled process is
\be\label{LdifftIRW2}
L^{\text{rw}}_N f(v) = \Delta_N f (v)=N^2\(f(v+ \tfrac{1}{N}) - 2 f(v) + f(v-\tfrac{1}{N})\), \qquad  v \in  \tfrac 1 N \Z
\ee
We denote by $(\caR_N,D(\caR_N))$ the Dirichlet form associated to the generator \eqref{LdifftIRW2}, that  is given by
\be\label{DirIRW}
\caR_N (f) = - \sum_{i \in \Z/N} f(i) \Delta_N f (i)  \mu_N(i)
\ee\label{disccountmun}
where $\mu_N$ is  the discrete counting measure on $\frac{1}{N} \Z$, this is
\be\label{munnnn}
 \mu_N(i) = \tfrac{1}{N}, \qquad \text{for all} \quad i \in \tfrac{1}{N} \Z
\ee
which is reversible for the dynamics. 
We are going to prove the Mosco convergence of the sequence of Dirichlet forms $\{(\caR_N,D(\caR_N))\}_N$ to  the Dirichlet form $(\caE_{\text{bm}},D(\caE_{\text{bm}}))$, i.e. the Dirichlet form associated to the standard Brownian motion in $\R$
\be
\caE_{\text{bm}}(f)= \frac{1}{2} \int_{\R} f\myprime(x)^2 dx.
\ee

\subsection*{Proof of Mosco convergence for RW}

\subsubsection*{Convergence of Hilbert spaces}
For the sequence of Hilbert spaces
\be\label{HrwN}
H^{\text{rw}}_N :=L^2 (\tfrac{1}{N}\Z, \mu_N )
\ee
where $\mu_N$ is as in \eqref{munnnn}. It is easy to see that we can guarantee the  convergence of $\{ H^{\text{rw}}_N \}_{N \geq 1}$ to the Hilbert space
\be\label{Hrw}
H^{\text{bm}}:=L^2 (\R,dx)
\ee
 i.e. the space of Lebesgue square-integrable functions in $\R$, by means  of the  restriction operators
\be\label{PHINIRW}
\{\Phi_N:C_k^{\infty}(\R) \subset H^{\text{bm}} \to H^{\text{rw}}_N\}_N \qquad \text{defined by}\qquad    \Phi_N f = f \mid_{\Z/N}.
\ee

\br
The choice of the space of all compactly supported smooth functions $C:=C_k^{\infty}(\R)$ as dense set for our Hilbert space turns out to be particularly convenient since it is a core  of the Dirichlet form associated to the Brownian motion. As a consequence,
we can make use of  the same set  also for proving that \eqref{equasecondass} is satisfied.
\er
\noindent
\subsubsection*{RW: Mosco I}

In order to prove that Assumption  \ref{strongcontinuity} is satisfied, it is convenient to split the proof in two cases depending whether $f$ belongs or not to the effective domain of $(-\Delta)^{-1/2}$. Hence, since $\Phi_N f\in H^{\text{rw}}_N$ is strongly convergent to $f \in C_k^\infty (\R)$, it is sufficient to prove Propositions \ref{StrongConPrepRW} and \ref{StrongConPrepRWinf} below:
\bp\label{StrongConPrepRW}
For any $f \in D((-\Delta)^{-1/2})$ we have
\be
\lim_{N \to \infty} \caR_N^* ( \Phi_N f ) = \caE_{{bm}}^* (f) \nn
\ee
\ep

\bpr
Let $G(x)$ be the Green's function of the Laplacian in $\R$, i.e. the fundamental solution to the problem
$\Delta G =\delta_0$
that is  given by
$G(x) = -\lvert x \rvert$.
We refer the reader to \cite{folland1995introduction} for more details on Green's functions.
Let $f$ be as in the statement, then,  by standard variational arguments  we know that
\beq
\caE_{\text{bm}}^* (f)  &=& \sup_{g \in D((-\Delta)^{1/2})} \left( \langle g,f \rangle  -\frac{1}{2}\norm{(-\Delta)^{1/2} g}_{L^2(\R)}^2   \right) = \frac{1}{4} \norm{ (-\Delta)^{-1/2} f }_{L^2(\R)}^2 \nn \\
&=& -\frac{1}{4} \langle f, G*f \rangle_{L^2(\R)}=  \frac{1}{4} \int_{\R} \int_{\R} f (x) f(y) \lvert x-y\rvert dx dy \label{E*} \nn
\eeq
Analogously, for the discrete case, we can write
\beq
\caR_N^* (\Phi_N f) = -\frac{1}{4}  \langle \Phi_N f, \Delta_N^{-1} \Phi_N f  \rangle_{H^{\text{rw}}_N} = -\frac{1}{4N}  \sum_{ i,j \in  \Z/N}  \Phi_N f(i) \cdot \Phi_N f(j) \cdot G_N(i-j) \nn \\
= -\frac{1}{4N}  \sum_{ i,j \in  \Z}  \Phi_N f(\tfrac{i}{N}) \cdot \Phi_N f(\tfrac{j}{N}) \cdot G_N(\tfrac{i-j}{N}) \nn
\eeq
where $G_N(\cdot)$ is the Green's function of the discrete Laplacian $\Delta_N$  in $\tfrac{1}{N}\Z$, i.e.  the  solution of the discrete problem:
\be\label{GN}
\Delta_N G_N = \delta_{0}\qquad \text{in }\Z/N
\ee
we refer to Chapter 5 in \cite{lawler2010random} for more details on discrete Green's functions. Notice that
\be
\tfrac 1 {N^2}\,{G}_1 (i) = G_N(\tfrac{i}{N})\qquad \forall i \in \Z \nn
\ee
where ${G}_1(\cdot)$ is the solution of \eqref{GN} for $N=1$,
then we can re-write
\beq\label{Greensapprox}
\caR_N^* (\Phi_N f)  = -\frac{1}{4N^3}  \sum_{ i,j \in  \Z} \Phi_N f(\tfrac{i}{N})\cdot \Phi_N f(\tfrac{j}{N}) \cdot {G}_1(i-j)
\eeq
By Theorem 4.48 in \cite{lawler2010random} we have that, for $i \neq j$, there exists $C, \beta >0$ such that
\be
{G}_1(i-j) = -\lvert i-j \rvert + C + O( e^{-\beta \lvert i-j \rvert}). \nn
\ee
Incorporating the above expression in \eqref{Greensapprox} we obtain
\beq\label{Greensapprox2}
\caR_N^* (\Phi_N f) = \tfrac{1}{4N^3}  \sum_{ \substack{i, j \in  \Z \\ i\neq j}} \Phi_N f(\tfrac{i}{N})  \Phi_N f(\tfrac{j}{N}) \( \lvert i-j \rvert + C + O( e^{-\beta \lvert i-j \rvert}) \) - \tfrac{1}{4N^3}  \sum_{ i \in  \Z}\(\Phi_N f(\tfrac{i}{N})\)^2 {G}_1(0) \nn
\eeq
notice that the sum on the diagonal vanishes as $N\to \infty$. Even more, thanks to the factor $N^{-3}$ in front of the two dimensional sum, we have that
\be
\lim_{N \to \infty} \frac{1}{4N^3}  \sum_{ \substack{i, j \in  \Z \\ i\neq j}} \Phi_N f(\tfrac{i}{N}) \Phi_N f(\tfrac{j}{N}) \(  C + O( e^{-\beta \lvert i-j \rvert}) \) = 0 \nn
\ee
Then we have
\beq\label{Greensapprox2}
&&\lim_{N \to \infty} \caR_N^* (\Phi_N f)  =  \lim_{N \to \infty} \frac{1}{4N^3}  \sum_{ \substack{i, j \in  \Z \\ i\neq j}} \Phi_N f(\tfrac{i}{N}) \cdot \Phi_N f(\tfrac{j}{N}) \cdot \lvert i-j \rvert   \\
&&=  \lim_{N \to \infty} \frac{1}{4N^2}  \sum_{ \substack{i, j \in  \Z \\ i\neq j}} \Phi_N f(\tfrac{i}{N}) \cdot \Phi_N f(\tfrac{j}{N})  \cdot \lvert\tfrac{i-j}{N} \rvert  = \frac{1}{4} \int_{\R} \int_{\R} f (x) f(y) \lvert x-y\rvert dx dy = \caE_{\text{bm}}^* (f)\nn
\eeq
which completes the proof.
\epr

\noindent
In order to conclude Assumption \ref{strongcontinuity} it remains to consider $f$ such that it does not belong to the domain of $ D((-\Delta)^{-1/2})$, this is $f$ such that 
$\caE_{\text{bm}}^* (f) = \infty$.

\bp\label{StrongConPrepRWinf}
For any $f \in H^{\text{bm}} \setminus D((-\Delta)^{-1/2})$ we have
$\lim_{N \to \infty} \caR_N^* ( \Phi_N f ) =  \infty$.
\ep

\bpr
Let $f $ be as in the statement, on the one hand we know
\beq
\caE_{{bm}}^* (f) &&= \frac{1}{4} \langle \Delta^{-1/2}f, \Delta^{-1/2}f  \rangle_{H^{\text{bm}}} = \frac{1}{4} \| \Delta^{-1/2}f   \|_{H^{\text{bm}}}^2 \nn \\
&&= \frac{1}{ 8 \pi} \| \widehat{\Delta^{-1/2}f }  \|_{L^2(dq)}^2 = \frac{1}{ 8 \pi} \| (iq)^{-1}\widehat{f}  \|_{L^2(dq)}^2 = \frac{1}{ 8 \pi} \int_{\R} \frac{(\widehat{f} (q))^2}{q^2} \; dq
\eeq
where $\widehat{f}$ denotes the Fourier transform of $f$. In the fourth line we used Plancherel's theorem, and in the fifth the differentiation property of the transform.\\

\noindent
Analogously for the discrete setting we have:
\beq\label{eqdualinfi1}
\caR_N^* ( \Phi_N f ) = \frac{1}{4} \langle \Phi_N f, \Delta_N^{-1} \Phi_N f  \rangle_{H_N^{\text{rw}}} = \frac{1}{4N^3}\sum_{ x \in \Z} \Phi_N f(\tfrac{x}{N})\cdot \Delta_1^{-1} \Theta_N f(\tfrac{x}{N})
\eeq
where $\Theta_N f: \Z \to \R$ is given by $\Theta_N f(x) := \Phi_N f(\tfrac{x}{N})$.
Let us denote by $\{  X_{t} : t \geq 0 \}$ the continuous time random walk on $\Z$ started at $ x$. Then we have that $\Delta_1^{-1} \Phi_N f(\frac{x}{N})$ is given by
\be\label{eqdualinfi2}
\Delta_1^{-1} \Phi_N f(\tfrac{x}{N}) = \frac{1}{2 \pi} \int_{-\pi}^{\pi} \frac{ \hat{f}_N(k) e^{-i k x} }{ 2 - 2 \cos k} dk
\ee
where
\beq
\widehat{f}_N(k)= \sum_{x \in \Z } \Theta_N f(x) e^{ ik x}
= \sum_{x \in \Z } \Phi_N f(\tfrac{x}{N}) e^{ ik x} \nn
\eeq
Substitution of \eqref{eqdualinfi2} in \eqref{eqdualinfi1} gives:
\beq
\caR_N^* ( \Phi_N f ) &&= \frac{1}{8 \pi N^3} \int_{-\pi}^{\pi} \frac{ \hat{f}_N(k) }{ 2 - 2 \cos k}  \sum_{ x \in \Z} \Phi_N f(\tfrac{x}{N}) e^{-i k x}\, dk \nn \\
&&= \frac{1}{ 8 \pi N^3} \int_{-\pi}^{\pi} \frac{ (\hat{f}_N(k))^2 }{ 2 - 2 \cos k}\, dk 
= \frac{1}{8 \pi} \int_{-\pi N}^{\pi N} \frac{ (\tfrac{1}{N}\hat{f}_N(\frac{q}{N}))^2 }{ N^2(2 - 2 \cos \tfrac{q}{N})}\, dq \nn
\eeq
at this point, in order to get convergence to the limitng dual we use the limits
\be
\lim_{N \to \infty} N^2 (2 - 2 \cos \tfrac{q}{N} ) =  q^2 \qquad \text{and} \qquad \lim_{ N \to \infty} \tfrac{1}{N} \widehat{f}_N(\tfrac{q}{N}) = \widehat{f}(q)
\ee
and 
by Fatou's Lemma we finish the proof.
\epr

\subsubsection*{RW: Mosco II}
For what concerns the second condition of Mosco convergence, we  choose $K:= C_k^\infty (\R)$ that is a core of $\caE_{\text{bm}}$. In this way, for all $f \in C_k^\infty (\R)$, we can consider the restrictions $\Phi_N f$ (strongly-convergent to $f$) and Taylor expand them to prove  that:
\beq
\lim_{N \to \infty} \caR_N (\Phi_N f)
&&=- \frac 1 N\lim_{N \to \infty}   \sum_{i \in \Z/N } \Phi_N f(i)  \Delta_N \Phi_N f (i)  \nn \\
&&=- \lim_{N \to \infty}  \frac{1}{N} \sum_{i \in \Z/N }  f(i) \Delta_N  f (i)  =- \lim_{N \to \infty}  \frac{1}{2N} \sum_{i \in \Z}  f(\tfrac{i}{N}) f\mydprime(\tfrac{i}{N}) + O(\frac{1}{N})  \nn \\
&&= -\frac{1}{2}\int_{\R} f(x) f\mydprime(x) dx = \frac{1}{2}\int_{\R} f\myprime(x)^2 dx =\caE_{\text{bm}}(f)
\eeq
which concludes the proof of Assumption 2. \qed

\br\label{extdMoscoRWBM}
Notice that Theorem 4.48 in \cite{lawler2010random} also applies for the finite range case and hence the results concerning Mosco convergence to the corresponding Brownian motion can be extended to the finite range setting modulus a multiplicative constant depending on the second moment of the transition $p$.
\er

\section*{Acknowledgements}
The authors would like to thank Mark Peletier for helpful discussions; The authors also would like to thank valuable comments from an anonymous reviewer.  M. Ayala acknowledges financial support from the Mexican Council on Science and Technology (CONACYT)  via the scholarship 457347.

\bibliographystyle{elsarticle-num}
\bibliography{Biblio}

\begin{thebibliography}{10}
\expandafter\ifx\csname url\endcsname\relax
  \def\url#1{\texttt{#1}}\fi
\expandafter\ifx\csname urlprefix\endcsname\relax\def\urlprefix{URL }\fi
\expandafter\ifx\csname href\endcsname\relax
  \def\href#1#2{#2} \def\path#1{#1}\fi

\bibitem{carinci2017exact}
G.~Carinci, C.~Giardina, F.~Redig, Exact formulas for two interacting particles
  and applications in particle systems with duality, arXiv preprint
  arXiv:1711.11283 (2017).

\bibitem{amir1991sticky}
M.~Amir, Sticky {Brownian} motion as the strong limit of a sequence of random
  walks, Stochastic processes and their applications 39~(2) (1991) 221--237.

\bibitem{cao2014dynamics}
J.~Cao, P.~Chleboun, S.~Grosskinsky, Dynamics of condensation in the totally
  asymmetric inclusion process, Journal of Statistical Physics 155~(3) (2014)
  523--543.

\bibitem{chau2015explosive}
Y.-X. Chau, C.~Connaughton, S.~Grosskinsky, Explosive condensation in symmetric
  mass transport models, Journal of Statistical Mechanics: Theory and
  Experiment 2015~(11) (2015) P11031.

\bibitem{beltran2017martingale}
J.~Beltr{\'a}n, M.~Jara, C.~Landim, A martingale problem for an absorbed
  diffusion: the nucleation phase of condensing zero range processes,
  Probability Theory and Related Fields 169~(3-4) (2017) 1169--1220.

\bibitem{grosskinsky2013dynamics}
S.~Grosskinsky, F.~Redig, K.~Vafayi, et~al., Dynamics of condensation in the
  symmetric inclusion process, Electronic Journal of Probability 18 (2013).

\bibitem{jatuviriyapornchai2020structure}
W.~Jatuviriyapornchai, P.~Chleboun, S.~Grosskinsky, Structure of the condensed
  phase in the inclusion process, Journal of Statistical Physics 178~(3) (2020)
  682--710.

\bibitem{mosco1994composite}
U.~Mosco, Composite media and asymptotic {Dirichlet} forms, Journal of
  Functional Analysis 123~(2) (1994) 368--421.

\bibitem{kuwae2003convergence}
K.~Kuwae, T.~Shioya, Convergence of spectral structures: a functional analytic
  theory and its applications to spectral geometry, Communications in analysis
  and geometry 11~(4) (2003) 599--674.

\bibitem{opoku2015coupling}
A.~Opoku, F.~Redig, Coupling and hydrodynamic limit for the inclusion process,
  Journal of Statistical Physics 160~(3) (2015) 532--547.

\bibitem{kipnis2013scaling}
C.~Kipnis, C.~Landim, Scaling limits of interacting particle systems, Vol. 320,
  Springer Science \& Business Media, 2013.

\bibitem{grosskinsky2011condensation}
S.~Grosskinsky, F.~Redig, K.~Vafayi, Condensation in the inclusion process and
  related models, Journal of Statistical Physics 142~(5) (2011) 952--974.

\bibitem{Grosskinsky2003CondensationIT}
S.~Grosskinsky, G.~M. Sch{\"u}tz, H.~Spohn, Condensation in the zero range
  process: Stationary and dynamical properties, Journal of Statistical Physics
  113 (2003) 389--410.

\bibitem{evans2005nonequilibrium}
M.~R. Evans, T.~Hanney, Nonequilibrium statistical mechanics of the zero-range
  process and related models, Journal of Physics A: Mathematical and General
  38~(19) (2005) R195.

\bibitem{fukushima1980dirichlet}
M.~Fukushima, Dirichlet forms and Markov processes, North-Holland Publishing
  Company, 1980.

\bibitem{kolesnikov2006mosco}
A.~V. Kolesnikov, Mosco convergence of {Dirichlet} forms in infinite dimensions
  with changing reference measures, Journal of Functional Analysis 230~(2)
  (2006) 382--418.

\bibitem{andres2010particle}
S.~Andres, M.-K. von Renesse, Particle approximation of the wasserstein
  diffusion, Journal of Functional Analysis 258~(11) (2010) 3879--3905.

\bibitem{howitt2007stochastic}
C.~J. Howitt, Stochastic flows and sticky brownian motion, Ph.D. thesis,
  University of Warwick (2007).

\bibitem{hewitt1975real}
E.~Hewitt, K.~Stromberg, Real and abstract analysis: a modern treatment of the
  theory of functions of a real variable (1975).

\bibitem{chen2012symmetric}
Z.-Q. Chen, M.~Fukushima, Symmetric Markov Processes, Time Change, and Boundary
  Theory (LMS-35), Vol.~35, Princeton University Press, 2012.

\bibitem{fukushima2010dirichlet}
M.~Fukushima, Y.~Oshima, M.~Takeda, Dirichlet forms and symmetric Markov
  processes, Vol.~19, Walter de Gruyter, 2011.

\bibitem{kato2013perturbation}
T.~Kato, Perturbation theory for linear operators, Vol. 132, Springer Science
  \& Business Media, 2013.

\bibitem{borodin2012handbook}
A.~N. Borodin, P.~Salminen, Handbook of Brownian motion-facts and formulae,
  Birkh{\"a}user, 2012.

\bibitem{folland1995introduction}
G.~B. Folland, Introduction to partial differential equations, Princeton
  university press, 1995.

\bibitem{lawler2010random}
G.~F. Lawler, V.~Limic, Random walk: a modern introduction, Vol. 123, Cambridge
  University Press, 2010.

\end{thebibliography}
\end{document}